%Variant February 21, 1997 
%By Viacheslav V. Nikulin. Steklov Math. Institute Moscow
%After February 1997 slava@nikulin.mian.su 

\input amstex
\documentstyle{amsppt}
\magnification=1200
\catcode`\@=11
\redefine\logo@{}
\catcode`\@=13

\define \bn{\Bbb N}
\define \bz{\Bbb Z}
\define \bq{\Bbb Q}
\define \br{\Bbb R}
\define \bc{\Bbb C}

\define \bk{\Bbb K}
\define \bo{\Bbb O}

\define \M{{\Cal M}}
\define\Ha{{\Cal H}}
\define\La{{\Cal L}}  
 
\define\0o{{\overline 0}}
\define\1o{{\overline 1}}
\define\rk{\text{rk}~}

%additional to the file 

\define\pr{\text{pr}}
\define\pd#1#2{\dfrac{\partial#1}{\partial#2}}

\define\NEF{\text{NEF}}
\define\Aut{\text{Aut}}
\define\arc{\text{arc}}

\TagsOnRight
%\NoBlackBoxes

\document

\topmatter
\title
K3 surfaces with interesting groups of automorphisms
\endtitle

\author
Viacheslav V. Nikulin \footnote{Supported by
Grant of Russian Fund of Fundamental Research and RIMS of Kyoto University 
\hfill\hfill}
\endauthor

\address
Steklov Mathematical Institute,
ul. Vavilova 42, Moscow 117966, GSP-1, Russia
\endaddress
\email
slava\@nikulin.mian.su
\endemail

\abstract
By the fundamental result of I.I. Piatetsky-Shapiro and I.R. Shafarevich 
(1971), the automorphism group $Aut(X)$ of a K3 surface $X$ over $\bc$ 
and its action on the 
Picard lattice $S_X$ are prescribed by the Picard lattice $S_X$.  
We use this result and our method 
(1980) to show finiteness of the set of Picard 
lattices $S_X$ of rank $\ge 3$ such that 
the automorphism group $Aut(X)$ of the K3 surface 
$X$ has a non-trivial invariant sublattice $S_0$ in $S_X$ where 
the group $Aut(X)$ acts as a finite group.   
For hyperbolic and parabolic lattices 
$S_0$ it has been proved by the author before (1980, 1995). 
Thus we extend this results to negative sublattices $S_0$.  

We give several examples of Picard lattices $S_X$ with parabolic and 
negative $S_0$.  
We also formulate the corresponding finiteness result for reflective 
hyperbolic lattices of hyperbolic type over purely real algebraic 
number fields.  

These results are important for the theory of Lorentzian Kac--Moody 
algebras and Mirror Symmetry.   
\endabstract

\rightheadtext
{Groups of automorphisms of K3 surfaces}
\leftheadtext{V.V. Nikulin}
\endtopmatter

\head 
0. Introduction
\endhead 

Here we want to extend theory of $2$-reflective and reflective 
hyperbolic lattices to so called hyperbolic type. For elliptic 
type this theory was developed by the author 
\cite{N2}---\cite{N14} and 
\'E. B. Vinberg (see e.g. \cite{V1}---\cite{V3}). For parabolic type 
it was developed in \cite{N11}. 

This theory is important for Lorentzian Kac--Moody algebras and 
Mirror symmetry. See \cite{B1}---\cite{B6}, \cite{GH}, \cite{GN1}---
\cite{GN7}, \cite{N1}---\cite{N14}. 

This paper was written during my stay at RIMS Kyoto University 
on 1996---1997. I am grateful to the Mathematical Institute for 
hospitality. 

\smallpagebreak 

A variant of this paper was published as preprint \cite{N15}.  

\smallpagebreak 

\head 
1. Automorphism groups of K3 surfaces and 2-reflective lattices of 
hyperbolic type  
\endhead

First we recall some general results about automorphism groups of K3 
surfaces. All of them one can find in the fundamental work by  
I.I. Piatetsky-Shapiro and I.R. Shafarevich  
\cite{P-SSh} where, in particular, automorphism groups of K3 surfaces 
over algebraically closed fields of characteristic $0$  
were described using the group of $2$-reflections of their Picard lattices.   

Let $X$ be a projective algebraic K3 
surface over an algebraically closed field $k$.  
Remind that it means: $X$ is a non-singular projective algebraic 
surface with the canonical class $K_X=0$ and 
$H^1(X, {\Cal O}_X)=0$. Let $S=S_X$ be the Picard lattice of $X$. 
The lattice $S$ is hyperbolic (has signature $(1,k)$).  
Let $V^+(S)$ be the half-cone of the light cone 
$V(S)=\{x \in S\otimes \br\,| \,x^2>0\}$ of $S$ which 
contains a class of hyperplane section (i.e. a polarization) of $X$ 
and $\La (S)=V^+(S)/\br_{++}$ the corresponding hyperbolic space with 
the distance $\cosh{\rho(\br_{++}x, \br_{++}y)}=(x,y)/\sqrt{x^2y^2}$, 
$x,y \in V^+(S)$. (Then the curvature is equal to $-1$.)   
The automorphism group 
$O^+(S)=\{\phi \in O(S)\,|\,\phi(V^+(S))=V^+(S)\}$ of the lattice $S$ 
is discrete in $\La (S)$ and has a fundamental domain 
of finite volume. We denote by $W^{(2)}(S)$ its subgroup generated 
by reflections in all elements with square $-2$ of $S$. Let 
$\NEF(X)$ be the nef cone of $X$.  
The convex polyhedron  $\M=\NEF(X)/\br_{++}\subset \La (S)$ is a 
fundamental chamber for the group $W^{(2)}(S)$, and the set 
$P(\M )$ of orthogonal vectors with square $-2$ to $\M$ coincides 
with classes of all exceptional curves on $X$ (all of them are  
irreducible non-singular rational). We have the semi-direct product 
$O^+(S)= W^{(2)}(S)\rtimes A(\M)$ where  
$A(\M)=\{\phi \in O^+(S)\,|\,\phi(\M)=\M\}$ is the group of symmetries 
of $\M$. The group $O(S)$ is an arithmetic group. 
From the general theory of arithmetic groups 
(e.g., see \cite{R}), it follows that if $\rk S\ge 3$, then
one can find a fundamental domain $\M_0\subset \M$ for $A(\M)$ 
acting in $\M$ which has {\it geometrically finite  volume}   
in $\La (S)$ (it means that it is a convex envelope of a finite set of 
points in $\La (S)$, possibly with some of them at infinity).  

The group $A(\M )$ is extremely important in 
geometry of the nef cone $\NEF(X)$ (and its dual Mori cone) 
of a K3 surface and geometry of the  
K3 surface $X$ itself. For example, it is known that 
the kernel of the natural homomorphism 
$\pi: \Aut (X)\to A(\M)$ is finite. Thus, if 
one has an information about the group $A(\M)$, he gets some 
information about the group $\Aut (X)$.   
For example, if the group $A(\M)$ is small, then the group $\Aut(X)$ 
is also small.  

One of the fundamental results of 
\cite{P-SSh} is that over the basic field $k$ of characteristic $0$ 
the homomorphism $\pi$ has also finite cokernel. 
In other words, the groups $\Aut (X)$ and $A(\M)$ 
coincide up to finite groups. Thus, in 
many cases one can replace the automorphism group $\Aut (X)$ 
by the group $A(\M)$. More exact results about $\pi$ one can 
find in \cite{N1}, \cite{N4}. We mention some of them. 
The kernel $G_0$ of $\pi$ is a cyclic group $\bz /n\bz$. If a 
prime $p|n$, then the lattice $S$ is $p$-elementary: $S^\ast/S$ is 
$(\bz/p\bz)^a$. In particular, if $n$ is divided by two different 
primes, the lattice $S$ is unimodular. For a general K3 surface $X$ 
with the Picard lattice $S$, the image of $\pi$ is equal to 
$A(\M)_0=\{g \in A(\M)\ |\ g|_{(S^\ast/S)}=\text{id}\}$. 
See more exact statements in \cite{N4, \S 10}.   

All questions we consider below don't change if one replaces 
the group $A(\M)$ by its subgroup of finite index. Below 
we work with the group $A(\M)$, but instead one can always consider 
the automorphism group $Aut(X)$ and its action on the Picard lattice 
$S$.   

Further we work with an arbitrary integral hyperbolic lattice $S$ and 
its group of $2$-reflections $W^{(2)}(S)$ and the fundamental polyhedron 
$\M$ for $W^{(2)}(S)$ acting in $\La(S)$.    
We leave to a reader a ``difficult problem'' of reformulation 
of our results below replacing the group $A(\M)$ by the automorphism 
group $Aut(X)$ in cases when the hyperbolic lattice $S$ is isomorphic to 
the Picard lattice of a K3 surface $X$. Over the field $k$ of 
characteristic $0$ a hyperbolic lattice $S$ is isomorphic to a Picard 
lattice of a K3 surface if and only if there exists 
a primitive embedding $S\subset L_{K3}$ where 
$L_{K3}$ is the even unimodular lattices of signature $(3,19)$ 
(in general, it follows from 
epimorphisity of the Torelli map proved by Vic.S. Kulikov 
in \cite{Ku}). In \cite{N2} simple necessary and sufficient 
conditions on $S$ are given to have this primitive embedding. For example,  
any even hyperbolic lattice $S$ of rank $\rk S\le 11$ 
has a primitive embedding to $L_{K3}$.

\definition{Definition 1} A hyperbolic lattice $S$ 
is called {\it $2$-reflective} if there exists 
a subgroup $G\subset A(\M)$ of 
finite index and an $G$-invariant sublattice  
$S_0\subset S$ such that $S_0\not=\{0\}$, 
and the group $G|_{S_0}$ is finite.  
There are three cases: 

1. The group $A(\M )$ has {\it elliptic type} if the lattice $S_0$ is 
hyperbolic (has one positive square). 
Obviously, this is equivalent for the group $A(\M )$ to be finite. 
Then the maximal $S_0$, one can take, is equal to $S_0=S$.  

2. The group $A(\M )$ has {\it parabolic type} if $A(\M )$ is 
infinite and $S_0$ has a one dimensional kernel 
$\bz c$ where $c^2=0$ and $c\not=0$. 
One can prove that then the element $c$ is fixed by the full 
group $A(\M )$ and $\br_{++}c \in \M$ (see \cite{N13}). The 
full group $A(\M )$ is trivial on the sublattice $S_0=[c]$.   
The maximal sublattice $S_0$ is equal to the primitive sublattice in 
$(c)^\perp_S$ generated by $c$ and all elements with square $-2$ in 
$(c)^\perp_S$.  

3. The group $A(\M )$ has {\it hyperbolic type} if $A(\M )$ has 
not elliptic or parabolic type and the lattice $S_0$ is 
negative. Replacing $S_0$ by a sublattice generated by the finite 
set of negative sublattices $g(S_0)$, $g \in A(\M )$, 
we can always suppose that 
the full group $A(\M )$ is finite on a non-zero negative sublattice 
$S_0$. Obviously, there exists a maximal negative sublattice $S_0$ in 
$S$ with this property. 

According to the type of the group $A(\M)$, the lattice $S$ is called 
$2$-reflective of elliptic, parabolic or hyperbolic type respectively. 
\enddefinition 

We mention that in Definition 1 of a $2$-reflective hyperbolic lattice,  
replacing $S_0$ by $(S_0)^\perp_S$, one can always drop the 
condition that $G$ is finite on $S_0$. 

If $\rk S=1$, then $S$ is obviously $2$-reflective. If $\rk S=2$, then 
the $S$ is $2$-reflective if and only if $S$ has a non-zero 
element with square $-2$ or $0$.  
It was proved in \cite{N3}---\cite{N9} 
that number of $2$-reflective  
hyperbolic lattices $S$ of elliptic type 
is finite, and all of them were classified. For $\rk S=4$ it 
was done by \'E.B. Vinberg (see \cite{N8}). Recently in 
\cite{N11} finiteness of the set of reflective hyperbolic 
lattices of parabolic type also was proved (there we considered the   
general case of reflective (not only $2$-reflective) lattices). 
Here we want to extend these our results for 
hyperbolic type using the same method we had applied 
for elliptic and parabolic type. 
Moreover, we can prove this finiteness for more general 
class of lattices $S$.  

\definition{Definition 2} A hyperbolic lattice $S$ is 
$2$-reflective of {\it hyperbolic type} if there exists a 
subgroup $G\subset A(\M)$ of finite index and a non-zero 
real negative subspace $K\subset S\otimes \br$ such that 
$g(K)=K$ for any $g\in G$. Equivalently, the group $G$ 
preserves a proper subspace $\La_0\subset \La(S)$ which is 
orthogonal to $K$. 

It seems, Definition 2 is more general 
than Definition 1. Obviously, if $S$ is $2$-reflective of hyperbolic 
type in the sense of Definition 1, 
then $S$ is $2$-reflective of hyperbolic type by Definition 2.   
Any hyperbolic lattice $S$ of the rank $\rk S \le 2$ is $2$-reflective in 
the sense of Definition 2.  
\enddefinition 

We prove   

\proclaim{Theorem 1} For a fixed rank $\rk S\ge 3$ number of 
$2$-reflective hyperbolic lattices $S$ of the hyperbolic type in 
the sense of Definition 2 is finite.
\endproclaim

\demo{Proof} Theorem 1 follows from 
Lemma 1 below. The same Lemma is valid for elliptic and parabolic type 
with the better constant 14 instead of $10+4\sqrt{5}$.

\proclaim{Lemma 1} Let $S$ be a $2$-reflective 
hyperbolic lattice of hyperbolic type and $n=\rk S \ge 3$. 
Then there are $\{\delta_1,...,\delta_n \}\subset P(\M )$ such that 

1) $\rk [\delta_1,...,\delta_n]=n$ 
\newline
and

2)
$
-2< (\delta_i, \delta_j)\le 
10+4\sqrt{5} \approx 18.94427191,\  1 \le i<j \le n. 
$
\endproclaim

Theorem 1 follows from Lemma 1 since number of 
lattices $[\delta_1,...,\delta_n]$ generated by 
$\{\delta_1,...,\delta_n\}$ is obviously finite 
(because $\delta_i^2=-2$ for any $\delta_i \in P(\M)$) 
and $[\delta_1,...,\delta_n]\subset S\subset 
[\delta_1,...,\delta_n]^\ast$.  

\demo{Proof of Lemma 1} First, suppose that $n\ge 4$. Equivalently, 
$\dim \La (S)\ge 3$. We consider two cases. 

Case 1. Suppose that there exists a hyperplane $\Ha_\delta$, 
$\delta \in P(\M)$, of a face (of the highest dimension) of 
$\M$ which does not intersect $\La_0$ including infinite points 
of $\La_0$. 

Let $p:\La (S)\to \La_0$ be the orthogonal projection into  
$\La_0$. The set $p(\overline{\Ha_\delta })\subset \La_0$ is compact. 
We consider an open ball $D$ in $\La_0$ which contains 
this set, and the orthogonal cylinder $C_{\overline{D}}$ 
over $\overline{D}$. 
The intersection $C_{\overline{D}}\cap \M$ 
has finite volume in $\La (S)$, and $\M$ is geometrically finite 
in $C_{\overline{D}}$ (i.e. $\M \cap C_{\overline{D}}$ is 
equal to intersection of $C_{\overline{D}}$ with a convex polyhedron in 
$\La (S)$ of geometrically finite volume). 
It follows that the face $\gamma $ of $\M$ 
containing in $\Ha_\delta$ has geometrically finite volume 
(in $\Ha_\delta$). Since $\dim \La (S)\ge 3$, 
there exists a $2$-dimensional face 
$\gamma_2\subset \gamma$ of geometrically finite volume. 
By Geometrical Lemma 3.2.1 in \cite{N5}, the plane polygon 
$\gamma_2$ has a {\it narrow part} 
which is given by four consecutive vertices $A_1A_2A_3A_4$ 
of $\gamma_2$ ($A_1=A_4$ if 
$\gamma_2$ is a triangle) such that the 
distance $\rho ((A_1A_2), (A_3A_4))< \text{arc\ } \cosh (7)$.  
We denote by $(AB)$ the  
line containing two different points $A$, $B$.   
(In fact, below we will prove a similar and even more 
complicated statement when we consider $n=3$.)  
We consider the subset $\Delta \subset P(\M)$ 
which consists of all $\delta \in P(\M)$ 
such that the hyperplane $\Ha_\delta$ 
contains one of lines $(A_1A_2)$, $(A_2A_3)$ or 
$(A_3A_4)$. Obviously, then the distance 
$\rho (\Ha_\delta, \Ha_{\delta^\prime })\le 
\text{arc}~\cosh (7)$ for any 
$\delta, \delta^\prime \in \Delta$. It follows 
that $(\delta, \delta^\prime)=
2\cosh \rho (\Ha_\delta, \Ha_{\delta^\prime})<14$. 
The set $\Delta$ has both properties 
1) and 2) of Lemma 1. The property 1) will be 
satisfied since all hyperplanes $\Ha_\delta$, $\delta \in \Delta $, 
don't have a common point and are not orthogonal to a hyperplane in 
$\La (S)$. Deleting unnecessary elements from $\Delta$, 
we get Lemma 1 for this case.  

Case 2. Suppose that all hyperplanes $\Ha_\delta$, $\delta \in P(\M)$,  
intersect the subspace $\La_0$. The polyhedron $\M$ has geometrically  
finite volume in the orthogonal cylinder $C_K$ over any compact subset 
$K\subset \La_0$. It follows that the polyhedron $\M$ has a 
vertex $v$ which is not contained in $\La_0$ (see considerations in 
\cite{N13, Lemma 1.2.2} where it is also proved 
that $\M$ contains a face $\gamma \subset \M$ of 
geometrically finite volume which is not contained in $\La_0$ and 
$\dim \gamma \ge \dim \La - \dim \La_0$. Consider a line $l$ which 
contains the vertex $v$ and is orthogonal to $\La_0$. Let $P=l\cap \La_0$. 
Consider a compact ball $D\subset \La_0$ with the center 
$P$ and the orthogonal cylinder $C_D$ over 
$D$. Since $\M$ has geometrically finite volume in $C_D$, 
there exists a hyperplane $\Ha_\delta$, 
$\delta \in P(\M)$, which intersects the line $l$ in a point $B$ which is 
different from $A$ and $\Ha_\delta$ does not contain $A$.  
Consider any hyperplane $\Ha_e$, $e \in P(\M)$, which contains 
the vertex $v$. Since both hyperplanes 
$\Ha_e$ and $\Ha_\delta$ intersect the 
subspace $\La_0$, it follows that 
$\rho (\Ha_e, \Ha_\delta )\le \arc \cosh 3$. 
It follows $(e,\delta)=2\cosh \rho(\Ha_e, \Ha_\delta)\le 6$. 
We consider the subset $\Delta \subset P(\M)$ 
which consists of the $\delta$ and all $e \in P(\M)$ 
such that the hyperplane $\Ha_e$ contains $v$.  
Again the set $\Delta$ obviously has both properties 1) and 2).  
It finishes the proof for $n\ge 4$.

Suppose that $n=3$ and the group $G$ is infinite. 
The $\La (S)$ is a hyperbolic plane 
and $\La_0$ is a line. 

The orthogonal projection of $\M$ into $\overline{\La_0}$ is an 
interval $[A,B]$ preserving by the group 
$G$. If one of points $A$ or $B$ is finite, then the group $G$ is 
finite. It follows that $\La_0\subset \M$.  
Let $e$ be an orthogonal vector to $\La_0$ with $e^2=-2$. 
All lines $\Ha_\delta$, 
$\delta \in P(\M)$, of 
sides of $\M$ (or the corresponding $\delta \in P(\M)$) are divided on two 
types depending on in what half-plane bounded by the 
line $\La_0$ they are contained. It is defined by the sign $(e, \delta)$. 
Since all lines $\Ha_\delta$, $\delta \in P(\M)$, may intersect 
the line $\La_0$ only at infinity, it follows that $|(e, \delta)|\ge 2$ 
for any $\delta \in P(\M)$.    
There exists also a possibility (which may really happen) that 
the line $\La_0$ is one of lines $\Ha_\delta$, $\delta \in P(\M)$, 
defining a side of $\M$. Then we assume that $e \in P(\M)$.  

We consider the invariant of $\M$ 
$$
a=\max \{|(e,\delta)|\ |\ \delta \in P(\M )\}.
$$
Geometrically $a=2\cosh{\rho}$ where $\rho$ is the maximal distance of 
lines of sides of the polygon $\M$ from the line $\La_0$. 
Since the group $G$ has finite index in $A(\M)$ and the fundamental 
domain for $G$ in $\M$ is obviously equal to the orthogonal cylinder 
with the base over the fundamental interval of $G$ acting on the line 
$\La_0$, it follows that the maximum $a$ always exists. 

We fix $f_2 \in P(\M)$ such that $|(f_2, e)|=a$. 
Replacing $e$ by 
$-e$, we can suppose that $a=(e,f_2)>0$. 
Since the line $\Ha_{f_2}$ may intersect the line $\La_0$ only at infinity, 
it follows that $a=(e,f_2)\ge 2$. Suppose that $a=(e,f_2)=2$. Then 
all lines $\Ha_\delta$, $\delta \in P(\M)$, have $|(\delta,e)|\le 2$. 
It then follows that they all intersect the line 
$\La_0$ in one of its 
infinite points. Then  number of these lines is 
finite (it is not greater than 4), and the group $G$ is finite. 
Thus, the invariant $a>2$, and the line $\Ha_{f_2}$ 
does not intersect the line $\La_0$. Consider the lines $\Ha_{f_1}$ and 
$\Ha_{f_3}$, $f_1, f_3 \in P(\M)$, which define the neighboring 
sides to the side $\Ha_{f_2}$ of $\M$. 
Thus, $\Ha_{f_1}$, $\Ha_{f_2}$ and $\Ha_{f_3}$ are lines of three 
consecutive sides of the fundamental polygon $\M$. 
Lemma 2 below is the analog of the statement above we used for 
a plane polygon $\gamma_2$ of geometrically finite volume.  

\proclaim{Lemma 2} We have: $(f_1,f_2)\le 2$, $(f_2,f_3)\le 2$ and 
$$
-2< (f_1, f_3)\le 14+{64\over a^2-4} \,.
$$  
\endproclaim

\demo{Proof} Since lines $\Ha_{f_1}$ and $\Ha_{f_2}$ intersect 
each other, $(f_1,f_2)\le 2$. Similarly, $(f_2,f_3)\le 2$.  
We have $(f_1,f_3)=2 \cosh{\rho(\Ha_{f_1},\Ha_{f_3})}$. Thus, 
we should estimate the distance $\rho(\Ha_{f_1},\Ha_{f_3})$. 
Consider the terminals $P$ of 
the line $\Ha_{f_1}$ and $Q$ of the line $\Ha_{f_3}$ on 
the absolute (i.e. the infinity of $\La(S)$), such that the 
line $\Ha_{f_2}$ is between lines 
$(PQ)$ and $\La_0$.   
Clearly, $\rho ((PQ),\La_0)\ge (\rho (\Ha_{f_2},\La_0)$.  
Let $f$ be the orthogonal vector with square $f^2=-2$ 
to the line $(PQ)$, and $(f, e)>0$. 
Then $b=(f, e)\ge (f_2,e)=a$. Thus, it is 
sufficient to prove Lemma replacing the element $f_2$ by 
$f$ and $a$ by $b$ respectively. Let us 
consider the perpendiculars $PP_1$ and $QQ_1$ to the line $\La_0$ 
where $P_1,Q_1 \in \La_0$. 
Consider the line $(P^\prime P)$ which is obtained from 
the line $(PQ)$ by the symmetry with respect to the line $(PP_1)$. 
Similarly, consider the line $(QQ^\prime)$ which is obtained from the 
line $(PQ)$ by the symmetry with respect to the line $(QQ_1)$. 
The lines $(P^\prime P)$ and $(QQ^\prime )$ have the same distance 
from $\La_0$ as the line $(PQ)$. 
Thus, they are further from the line $\La_0$ than 
the lines $\Ha_{f_1}$ and $\Ha_{f_3}$ because $(f_1,e)\le a \le b$ and 
$(f_3,e)\le a \le b$. It then follows that any interval with terminals on 
the lines $(P^\prime P)$ and $QQ^\prime$ intersects both lines 
$\Ha_{f_1}$ and $\Ha_{f_3}$. Thus, 
$\rho((P^\prime P),(QQ^\prime))\ge \rho (\Ha_{f_1},\Ha_{f_3})$. Let 
$f_1^\prime$ and $f_3^\prime$ are orthogonal vectors with square 
$(f_1^\prime)^2=(f_3^\prime)^2=-2$ to the lines 
$(P^\prime P)$ and $(QQ^\prime )$ respectively, and $(f_1^\prime,e)>0$,  
$(f_3^\prime,e)>0$. Then $(f_1^\prime,f_3^\prime )\ge (f_1,f_3)$, and 
it is sufficient to estimate $(f_1^\prime, f_3^\prime)$. 
Let $g$ and $h$ are orthogonal vectors with square $-2$ 
to lines $PP_1$ and $QQ_1$ and directed outwards from the quadrangle 
$P_1PQQ_1$. Then the vectors $f_1^\prime$ and $f_3^\prime$ are 
given by the reflections in $g$ and $h$ respectively. 
We have $f_1^\prime=f+2g$ and $f_3^\prime =f+2h$. Thus,  
$(f_1^\prime,f_3^\prime)=f^2+2(f,h)+2(f,g)+4(g,h)=
6+4x$ where $x=(g,h)$. The Gram matrix of the four 
vectors $e,g,f,h$ is equal to 
$$
\Gamma=\left(\matrix
-2 & 0 & b & 0\\
 0 &-2 & 2 & x\\
 b & 2 & -2& 2\\
 0 & x &  2&-2
\endmatrix\right). 
$$
We have $\det(\Gamma )=x^2(b^2-4)-16x-4b^2-16=0$ and  
$x=2+16/(b^2-4)$. Thus 
$(f_1,f_3)\le (f_1^\prime,f_3^\prime)=14+64/(b^2-4)\le 14+64/(a^2-4)$. 
It proves Lemma 2.  
\enddemo  

We denote by $a_0=2\sqrt{2+\sqrt{5}}=4.116342....$ 
the root of the equation 
$$
x^4-16x^2-16=0.
$$

We consider three cases. 

Case 1. Suppose that $a\ge a_0$. 

Then elements 
$f_1,f_2,f_3\in P(\M)$ of Lemma 2 give elements for the 
Lemma 1 we are looking for.

Case 2. Suppose that $a \le a_0$ and the line $\La_0$ is one of lines 
of sides of $\M$. Equivalently, $e \in P(\M)$. 

Consider a vertex 
$v$ of $\M$. We take two lines $\Ha_{f_1}$ and $\Ha_{f_2}$, 
$f_1,f_2 \in P(\M)$, of the two consecutive sides of $\M$ which contain 
the vertex $v$. Then elements 
$f_1,f_2,f_3=e$ give elements of Lemma 1 we are looking for. Really, 
$-2< (f_1,f_2)\le 2$ because lines $\Ha_{f_1}$ and $\Ha_{f_2}$ intersect 
each other. We have $(f_1,f_3)\le a \le a_0\le 18$, 
$(f_2,f_3)\le a \le a_0\le 18$  by definition of the invariant $a$.     

Case 3. Suppose that $a \le a_0$ and the line $\La_0$ is not a line 
of a side of $\M$. 

Like for the case 2, we take elements 
$f_1, f_2 \in P(\M)$ such that $\Ha_{f_1}$ and $\Ha_{f_2}$ contain 
a vertex $v$ of the polyhedron $\M$. We can suppose that 
$(f_1,e)>0$ and $(f_2,e)>0$. Consider a perpendicular $(vP)$ 
to the line $\La_0$. The line $(vP)$ intersects another line $\Ha_{f_3}$, 
$f_3 \in P(\M)$, of a side of $\M$. This line is contained in the second 
half-plane bounded by $\La_0$. Thus, the line $\La_0$ is between lines 
$\Ha_{f_1}$, $\Ha_{f_3}$, and $\Ha_{f_2}$, $\Ha_{f_3}$ as well. 
Suppose that the line $(vP)$ is equal to $(RS)$ where $R$, $S$ 
are points on the absolute and $v$ is between $R$ and $P$. 
Consider lines $l_1$ and $l_2$ 
with orthogonal vectors with 
square $-2$ $w_1$ and $w_2$ respectively and such that 
$l_1$ contains $R$, 
$l_2$ contains $S$, lines $l_1$, $l_2$ are contained in different 
half-planes bounded by $(RS)$ and in different half-planes bounded by 
$\La_0$, and $(w_1,e)=a$, $(w_2,e)=-a$. Then $(w_1,w_2)\ge (f_1,f_3)$ and 
$(w_1,w_2)\ge (f_2,f_3)$ because $0< (f_1,e)\le a$, 
$0< (f_2,e)\le a$ and 
$-a \le (f_3,e)<0$ by definition of the invariant $a$. (Geometrically, 
lines $l_1$ and $l_2$ are further from each other than the lines 
$\Ha_{f_1}, \Ha_{f_3}$ and lines $\Ha_{f_2}$, $\Ha_{f_3}$ as well. One 
should draw the corresponding picture to make sure.) 
Lines $l_1$ and $l_2$ are symmetric with respect to 
the point of intersection of lines 
$\La_0$ and $(RS)$.  If $g$ is the orthogonal vector with 
$g^2=-2$ to the line $(RS)$ and $(g,w_1)=2$, 
then $w_2$ is obtained from $w_1$ by the composition of 
symmetries in $g$ and $e$. It follows that 
$w_2=w_1+2g+ae$ and $(w_1,w_2)=2+a^2$. Thus, we have proved that 
$(f_3,f_1)\le 2+a^2$ and $(f_3,f_2)\le 2+a^2$. 
Since $a\le a_0$,  
elements $f_1,f_2,f_3$ give elements of Lemma 1 we are looking for.    

The equation for $a_0$, we have introduced above, is the equation 
$2+x^2=14+64/(x^2-4)$ when Lemma 2 and considerations of Case 3 give 
the same result.  

It finishes the proof of Lemma 1 and Theorem 1. 
\enddemo 
\enddemo

\head 
2. Reflective hyperbolic lattices 
\endhead

Let $\bk$ be a purely real algebraic number field of finite degree 
$[\bk : \bq]$ and $\bo$ its ring of integers. A {\it lattice} 
$S$ over $\bk$ is a projective module over $\bo$ of a 
finite rank equipped with 
a non-degenerate symmetric bilinear form with values in $\bo$. 
A lattice $S$ is called {\it hyperbolic} if there exists an embedding 
$\sigma^{(+)}:\bk \to \br$ (it is called the {\it geometric embedding}) 
such that the real symmetric bilinear form 
$S\otimes \br$ is hyperbolic (i.e. it has exactly $1$ positive 
square), and for all other embeddings $\sigma \not=\sigma^{(+)}$ the 
real form $S\otimes \br$ is negative definite. Like for integral 
hyperbolic lattices, using the geometric embedding, 
one can define the hyperbolic space $\La (S)$, 
the groups $O^+(S)$, the group  
$W(S)$ generated by all reflections of the lattice $S$ 
(i.e. automorphisms from $O^+(S)$ 
which act as reflections with respect to hyperplanes in 
$\La (S)$), fundamental 
polyhedron $\M\subset \La (S)$ for $W(S)$, 
set $P(\M)\subset S$ of orthogonal vectors 
to faces of $\M$, and the group of symmetries $A(\M)$ of $\M$.  

Applying Definitions 1 and 2 to this situation, we get definition of 
{\it reflective hyperbolic lattices} of {\it elliptic, parabolic and 
hyperbolic types}. Here a sublattice $S_0\subset S$ (over $\bo$) is 
called hyperbolic, parabolic or negative if it generates a real subspace of 
this type for the geometric embedding. 
For parabolic type, one should consider $\bo$ instead 
of $\bz$; it then follows that $\bk=\bq$.  

Using combination of considerations here for the proof of Theorem 1 and 
considerations in \cite{N5} and \cite{N6}, we get 

\proclaim{Theorem 2} i) For a fixed $N$, the set of 
fields $\bk$ such that $[\bk:\bq]=N$ and 
there exists a hyperbolic reflective lattice 
$S$ over $\bk$ of rank $\rk S\ge 3$ 
(of any type: elliptic, parabolic or hyperbolic) is finite.  

ii) For a fixed rank $n \ge 3$ and a 
fixed purely real algebraic number field $\bk$, the set of 
reflective hyperbolic lattices (of any type: elliptic, parabolic 
or hyperbolic) $S$ over $\bk$ and of 
the rank $n$ is finite up to multiplication 
of the form of $S$ by $m \in \bk$. 
\endproclaim

For reflective hyperbolic lattices of elliptic type it was proved in 
\cite{N5} and \cite{N6}. For parabolic type it is proved in \cite{N11}. 
Thus, Theorem 2 generalizes that results for hyperbolic type. Like 
for elliptic and parabolic type, the proof is based on 

\proclaim{Lemma 3} For any reflective hyperbolic lattice $S$ of rank 
$\rk S=n\ge 3$  
there are $\delta_1,...,\delta_n\in P(\M)$ such that 

1) $\rk [\delta_1,...,\delta_n]=n$;

2) the Gram diagram of $\delta_1,...,\delta_n$ is connected (i.e. one 
cannot divide the set $\{\delta_1,...,\delta_n\}$ on two non-empty 
subsets which are orthogonal to one another);

3) for the geometric embedding,  
$$
{4(\delta_i,\delta_j)^2\over \delta_i^2\delta_j^2}< 200^2,\ 1\le i<j\le n.
$$
\endproclaim

\demo{Proof} For elliptic type Lemma 3 was proved in \cite{N5}, 
for parabolic type in  \cite{N11}  
with the constant $62$ instead of $200$. 
For hyperbolic type, the proof is similar and uses combination 
of considerations in \cite{N5} and considerations 
here for the proof of Lemmas 1 and 2 and Theorem 1. In fact, 
like in \cite{N13} and \cite{N11}, 
one can introduce a class of convex locally 
finite polyhedra $\M$ of {\it restricted hyperbolic type} and prove 
Lemma 3 for these polyhedra like we did it for 
polyhedra of elliptic and restricted parabolic type. 
Lemma 3 for this class of polyhedra 
is important as itself: for example, for Mori polyhedra 
(e.g. see \cite{N13} )  
and for Borcherds type automorphic products (see \cite{B5}, 
\cite{N11} and \cite{GN5}).   
We hope to present details of the proof somewhere.   
\enddemo 

We mention that Theorem 2 partly gives an answer to our 
question in \cite{N9, Sect.3}: to generalize results 
about crystallographic reflection groups in hyperbolic spaces 
for so called generalized crystallographic reflection groups in 
hyperbolic spaces.

\head
3. Some examples of $2$-reflective hyperbolic lattices 
\endhead

Here we present some examples of 2-reflective  
hyperbolic lattices over $\bz$.

\subhead 
3.1. $2$-elementary $2$-reflective hyperbolic lattices 
and K3 surfaces with involution 
\endsubhead
This example had been first considered in \cite{N4} (and, in fact, 
in \cite{N2}, \cite{N3} before, see also related results in 
\cite{N7}, \cite{N9}).  
Let $S$ be a $2$-elementary even 
hyperbolic lattice and $S\subset L_{K3}$. We remind that $S$ is 
called $2$-elementary if $S^\ast/S\cong (\bz /2\bz)^a$. 
By the Global Torelli Theorem \cite{P-SSh}, there exists a K3 surfaces 
$X$ over $\bc$ 
with $S=S_X$ and with the {\it canonical involution} $\sigma$ which 
acts trivially on $S_X$ and as multiplication by $-1$ on the transcendental 
lattice $T_X=S_X^\perp$. Obviously, the 
$\Aut(X)$ normalizes the involution $\sigma$. In particular, 
$\Aut (X)$ preserves the set $X^\sigma$ of points of $X$ 
fixed by this involution. The set 
$X^\sigma$ is a non-singular curve whose components 
generate the 
sublattice $S_0\subset S$ which is obviously invariant with respect to 
$Aut(X)$, and the group $Aut (X)$ is finite on $S_0$. Thus, if $X^\sigma 
\not=\emptyset$, the lattice $S$ is $2$-reflective. 

In \cite{N4} (actually, it had been done in \cite{N2} and \cite{N3}) 
the set $X^\sigma$ and the lattice $S_0$  were 
described explicitly using invariants of the 
lattice $S$, and all $2$-elementary lattices $S$ having a primitive 
embedding to $L_{K3}$ also were described. 
Let $r=\rk S$. For $S\cong U(2)\oplus E_8(2)$ the set $X^\sigma =\emptyset$ 
and the lattice $S$ is not 2-reflective (because it does not have 
elements with the square $-2$). (Here we denote by $K(t)$ a lattice $K$ 
with the form multiplied by $t \in \bq$. The lattice $U$ is the standard 
even unimodular lattice of signature $(1,1)$. The lattice $E_8$ is the 
standard even unimodular lattice of signature $(0,8)$.)    
For $S\cong U\oplus E_8(2)$, the set $X^\sigma=C_1^{(1)}\coprod C_1^{(2)}$ 
and the lattice $S$ is $2$-reflective of parabolic type.  
If $S$ is different from the two lattices above, then 
$X^\sigma =C_g\coprod C_0^{(1)}\coprod \dots \coprod C_0^{(k)}$ where 
$g=(22-r-a)/2$, $k=(r-a)/2$. Here the below index denotes the genus $g$ 
(or $0$)  
of the corresponding smooth irreducible curve, 
equivalently $C_g^2=2g-2$ and 
$(C_0^{(t)})^2=-2$, $1\le t \le k$.  Obviously, $g+k \le 11$ 
and excluding some exceptions for $g+k=10$ and $g+k=11$,  
all other cases $g,k$ satisfying 
$g\ge 0$, $k \ge 0$ and $g+k \le 11$ correspond to the K3 surfaces with 
involution.  

In particular, if $g>1$, the $S$ has elliptic type. If $g=1$, the lattice 
$S$ has elliptic or parabolic type. It was shown in \cite{N4} that for 
$g=1$ only the lattice 
$S=U\oplus E_8\oplus E_8\oplus \langle -2 \rangle$ has elliptic type. 
We remark that for 
$g=1$ the K3 surface $X$ has the canonical elliptic fibration $|C_1|$ 
which is preserved by the group $\Aut (X)$. 

If $g=0$, then $S$ 
is $2$-reflective of hyperbolic or parabolic type. 

Thus, we get a lot of examples of $2$-reflective hyperbolic 
lattices $S$ of elliptic, parabolic and  hyperbolic type of all 
possible rank $\rk S=1,\,2,\,\dots ,\,20$.   
All $2$-elementary hyperbolic lattices $S$ having a primitive 
embedding to $L_{K3}$ are $2$-reflective except the lattice 
$U(2)\oplus E_8(2)$. 

\smallpagebreak 

\subhead
3.2. A general result about $2$-reflective hyperbolic lattices $S$ 
of $\rk S \ge 6$ 
\endsubhead 

In \cite{N4} we proved that for any K3 surface $X$ with $\rk S\ge 6$ 
the $X$ always has an elliptic 
fibration with infinite automorphism group if the full automorphism group 
$Aut (X)$ is infinite.  It is interesting that this is true 
for arbitrary characteristic of the ground 
field. It then follows that 
for parabolic and hyperbolic types of Definition 1   
the canonical lattice $S_0$ is equal to intersection of all 
primitive sublattices in $S$ generated by  
components of fibers of elliptic fibrations on $X$ 
having an infinite automorphism group. (In \cite{N4}
we also claimed that $S_0$ is generated up to finite index 
by irreducible curves in $X$, but this statement is might be wrong.)    
If the lattice $S_0\not=0$, then $S$ is $2$-reflective of parabolic 
type if $S_0$ is parabolic, and $S$ is $2$-reflective of hyperbolic 
type if $S_0$ is negative. It shows 
geometrical sense of $S_0$ when $X$ has elliptic fibrations with 
infinite automorphism groups. The same statement is valid for all 
hyperbolic lattices 
$S$ of $\rk S\ge 6$ if one replaces the group $\Aut (X)$ 
by $A(\M)$ and elliptic fibrations by automorphism groups of infinite 
vertices of $\M$. 

We mention that using this idea one can prove that $2$-reflective 
(and reflective) hyperbolic lattices over $\bz$ of hyperbolic type 
do not exist in high dimension. This is known for elliptic type 
\cite{N3}, \cite{N4}, \cite{N6} and \cite{V2} and parabolic type 
\cite{N11}.

\subhead
3.3. $2$-reflective hyperbolic lattices of rank 3 
\endsubhead

We consider hyperbolic lattices $S_k=U\oplus \langle -2k \rangle$, 
$k \in \bn$. We consider 
the standard basis $e_1, e_2$ of the lattice $U$ and the 
standard bases $e_3$ of $\langle -2k \rangle$, i. e. with the Gram matrix 
$$ 
(e_i, e_j)=
\left(\matrix
0 & 1 & 0\\
1 & 0 & 0\\
0 & 0 & -2k
\endmatrix
\right).
$$
Consider $k=2$. 
For the fundamental polyhedron $\M_0$ of the full 
reflection group $W(S_2)$   
we have (using Vinberg's algorithm \cite{V1})  
$$
P(\M_0)_{\pr}=\{\delta_{0,1}=(0,0,1),\,\delta_{0,2}=(2,0,-1),\, 
\delta_1=(-1,1,0)\}
$$
with the Gram matrix 
$$
G(P(\M_0)_{\pr})= 
\left(\matrix
-4 & 4 & 0\\
4  &-4 & 2\\
0  & 2 & -2
\endmatrix\right).
$$
Here the first index $i$ in $\delta_{i,j}$ shows appearing of 
this element on the $i$-th step of the Vinberg's algorithm with 
the center $\rho=(1,0,0)$. Index $\pr$ shows that we consider 
primitive elements. We denote by $\Delta^{(2)}(S)$ the set of all 
elements of $S$ with square $-2$. 

All intermediate lattices $S_2 \supset S=S_{2,l} \supset [\Delta^{(2)}(S_2)]$ 
of index $l=[S_2:S]$ are equal to 
$$
\split 
S_2=S_{2,1}=[\delta_{0,1}, \delta_{0,2}, \delta_1,
(\delta_{0,1}+\delta_{0,2})/2] 
&\supset S_{2,2}=[\delta_{0,1}, \delta_{0,2}, \delta_1]\\
&\supset
S_{2,4}=[2\delta_{0,1}, \delta_{0,2}, \delta_1]=[\Delta^{(2)}(S_2)]
\endsplit, 
$$
and for all these lattices  
$$
P(\M_0)=\{\delta_{0,1}=(0,0,1),\,\delta_{0,2}=(2,0,-1),\, 
\delta_1=(-1,1,0)\}. 
$$
One should multiply that elements by some constants from $\bn$ to 
get elements from $S_{2,l}$. 
For all these lattices 
$$
P(\M)=A(\M )(\delta_1),\ \ 
A(\M )=[s_{\delta_{0,1}}, s_{\delta_{0,2}}]\cong D_\infty.
$$
We denote by $D_\infty$ the group on a line generated by two different 
reflections. Here $s_\delta$ denote the reflection in $\delta$. 
Thus, the group $A(\M)$ and the corresponding 
group $Aut(X)$ of a K3 surfaces $X$ with 
$S_X\cong S_{2,l}$ are isomorphic to $\bz$. 
It follows that all lattices $S_{2,l}$ are $2$-reflective of 
parabolic type because  
the group $A(\M)$ fixes the parabolic sublattice $\bz c$ where 
$c=\delta_{0,1}+\delta_{0,2}$ with $c^2=0$. 
The linear system $|c|$ defines the 
unique elliptic fibration on the K3 surface $X$. Moreover, this elliptic 
fibration has infinite automorphism group 
(it is $\bz$). Thus, we get three examples of 
K3 surfaces with the parabolic automorphism group $Aut(X)\cong \bz$ 
and with very nice geometry.   

One can check that 
$$
S_2\cong U \oplus \langle -4 \rangle;\  S_{2,2}\cong \langle 4 \rangle \oplus 
\langle -2 \rangle^2;\  
S_{2,4}\cong \langle 16 \rangle 
\oplus \langle -2 \rangle^2.
$$
 
We get similar examples for $k=3$, $5$, $7$, $13$. 
The case $k=3$ gives three 
lattices $S_3\cong U\oplus \langle -6 \rangle$; 
$S_{3,3}\cong \langle 18 \rangle \oplus A_2$; 
$S_{3,6}\cong \langle 72 \rangle \oplus A_2$. 
The case $k=5$ gives the only lattice $S_5\cong U\oplus \langle -10 \rangle$. 
The case $k=7$ gives two lattices $S_7=U\oplus \langle -14 \rangle$ 
and $S_{7,2}=[\Delta^{(2)}(S_7)]$. The case $k=13$ gives the only lattice 
$U\oplus \langle -26 \rangle$. Here one should use calculations in 
Sect. 4 below. 

Using the general method described above for the proof 
of Theorem 1, we can prove   

\proclaim{Theorem 3} Lattices $S_{k,l}$, 
$k=2$, $3$, $5$, $7$, $13$, defined above, 
are the only Picard lattices of rank $3$ of K3 surfaces $X$ 
such that $X$ has the only one elliptic fibration, moreover this 
fibration has an infinite automorphism group (it is then isomorphic to 
$\bz$).
\endproclaim 

\smallpagebreak 

Now we can look for all other $2$-reflective lattices in the series 
$S_k=U\oplus \langle -2k \rangle$ and their intermediate lattices 
$S_k\supset S_{k,l}\supset [\Delta^{(2)}(S_k)]$. We get the 
following and the only following $2$-reflective lattices: 
The lattice 
$S_1=U\oplus \langle -2 \rangle$ is $2$-reflective of  elliptic type. 
All other lattices will be $2$-reflective of parabolic type. They are: 
the lattice $S_{4}$ and its intermediate sublattices 
$S_{4,1}$, $S_{4,2}$, $S_{4,4}$, $S_{4,8}=[\Delta^{(2)}(S_4)]$. 
The lattice $S_{9}$ and its intermediate sublattices 
$S_{9,3}=[\Delta^{(2)}(S_9)]$.  
The lattice $S_{25}$. (Here one should use calculations in 
Sect. 4 below.) Lattices $S_4$, $S_9$ and $S_{25}$ and 
their intermediate sublattices give Picard lattices or rank $3$ 
of K3 surfaces having the unique elliptic 
fibration with infinite automorphism group and some other elliptic 
fibrations with finite automorphism group. One need to add very few 
cases to get the full list of K3 surfaces of this type. 

\smallpagebreak 

Consider the lattice $S=U(11)\oplus \langle -2 \rangle$ with 
the standard bases. 
The fundamental polyhedron $\M_0$ for $W(S)$ has 
the set of primitive orthogonal 
vectors $\delta_1,...,\delta_6$ with coordinates given by lines 
of the matrix $P(\M_0)$ and with the Gram matrix $G(P(\M_0))$ given 
below   
$$
P(\M_0)=\pmatrix{1}&{0}&{-1}\cr{0}&{0}&{1}\cr{-1}&{1}&{0}\cr
{3}&{3}&{-10}\cr{8}&{6}&{-23}\cr{5}&{2}&{-11}\cr\endpmatrix,
\ \ \ 
G(P(\M_0))=
\pmatrix{-2}&{2}&{11}&{13}&{20}&{0}\cr
{2}&{-2}&{0}&{20}&{46}&{22}\cr
{11}&{0}&{-22}&{0}&{22}&{33}\cr
{13}&{20}&{0}&{-2}&{2}&{11}\cr
{20}&{46}&{22}&{2}&{-2}&{0}\cr
{0}&{22}&{33}&{11}&{0}&{-22}\cr\endpmatrix . 
$$
It follows that the fundamental polyhedron $\M$ for $W^{(2)}(S)$ 
has 
$$
P(\M)=A(\M)(\{\delta_1,\delta_2,\delta_4, \delta_5\}),\ \  
A(\M)=[s_{\delta_3}, s_{\delta_6}]\cong D_\infty.
$$
The lattice $S$ is $2$-reflective of hyperbolic type since 
$A(\M)$ keeps the sublattice $[\delta_3,\delta_6]^\perp$ generated 
by $w=(2,2,-7)$ with $w^2=-10$. Thus, the lattice $S$ is 
$2$-reflective of hyperbolic type. A K3 surface $X$ over $\bc$ with 
the Picard lattice $S\cong U(11)\oplus \langle -2 \rangle$ has 
the automorphism group $\Aut (X)\cong \bz$ but this group is 
not the automorphism group of any elliptic fibration on $X$. All 
elliptic fibrations on $X$ have finite automorphism groups. We 
remark that the lattice $S=U(11)\oplus \langle -2 \rangle$ is 
$11$-dual to the lattice $S_{11}=U\oplus \langle -22 \rangle$ which 
we will consider in Sect. 4 below.   

Similarly one can check that the lattices $U(15)\oplus \langle -2 \rangle$ 
and $U(24)\oplus \langle -2 \rangle$ are $2$-reflective 
of hyperbolic type. K3 surfaces over $\bc$ with these Picard lattices 
have the group $\Aut (X)\cong \bz$ but this group is not a group of 
automorphisms of 
any elliptic fibration on $X$. All elliptic fibrations on $X$ have 
finite automorphism groups. We remark that $U(15)\oplus \langle -2 \rangle$ 
is $15$-dual to the lattice $S_{15}=U\oplus \langle -30 \rangle$, and 
the lattice $U(24)\oplus \langle -2 \rangle$ is dual to the lattice 
$S_{24}=U\oplus \langle -48 \rangle$ which we consider below.   

We hope to classify all $2$-reflective hyperbolic lattices 
of rank $3$ later. For elliptic type of rank 3 this had been done 
in \cite{N8}.

\head
4. Reflective hyperbolic lattices of rank 3 representing $0$  
\endhead

Here we consider the fundamental series $S_k=U\oplus \langle -2k \rangle$, 
$k \in \bn$, of hyperbolic lattices of rank $3$ representing $0$. 

Below we give results of our calculations using Vinberg's algorithm 
\cite{V1} of the fundamental polyhedron 
$\M$ for the full group $W(S_k)$ generated by reflections, if $k \le 60$.  
We use the standard bases $e_1$, $e_2$, $e_3$ 
introduced above. The first matrix (or the 
matrices $e$, $f$ for the hyperbolic type) gives the set $P(\M)$ of 
primitive orthogonal vectors to $\M$. The second matrix (or 
$G(e)$, $G(f)$ for the hyperbolic type) gives their Gram matrix. 
For the hyperbolic type the vector $w$ 
generates a negative sublattice $S_0$ in $S_k$ which is invariant 
with respect to  
the group $A(\M)$ which is given by the set of generators $C_i$. 
For hyperbolic type the set $P(\M)=A(\M)(e\cup f)$ or 
$P(\M)=A(\M)(e)$ (if we don't give the set $f$).

\noindent
$n=1: \hskip20pt 
\pmatrix{1}&{0}&{-1}\cr{0}&{0}&{1}\cr{-1}&{1}&{0}\cr\endpmatrix 
\ \ \ 
\pmatrix{-2}&{2}&{1}\cr{2}&{-2}&{0}\cr{1}&{0}&{-2}\cr\endpmatrix 
$. 

\noindent 
$n=2: \hskip20pt  
\pmatrix{2}&{0}&{-1}\cr{0}&{0}&{1}\cr{-1}&{1}&{0}\cr\endpmatrix 
\ \ \ 
\pmatrix{-4}&{4}&{2}\cr{4}&{-4}&{0}\cr{2}&{0}&{-2}\cr\endpmatrix  
$.

\noindent
$n=3: \hskip20pt 
\pmatrix{3}&{0}&{-1}\cr{0}&{0}&{1}\cr{-1}&{1}&{0}\cr\endpmatrix
\ \ \ 
\pmatrix{-6}&{6}&{3}\cr{6}&{-6}&{0}\cr{3}&{0}&{-2}\cr\endpmatrix 
$.

\noindent
$n=4: \hskip20pt 
\pmatrix{4}&{0}&{-1}\cr{0}&{0}&{1}\cr{-1}&{1}&{0}\cr\endpmatrix
\ \ \ 
\pmatrix{-8}&{8}&{4}\cr{8}&{-8}&{0}\cr{4}&{0}&{-2}\cr\endpmatrix$ . 

\noindent
$n=5: \hskip20pt  
\pmatrix{5}&{0}&{-1}\cr{0}&{0}&{1}\cr{-1}&{1}&{0}\cr{2}&{2}&{-1}\cr
\endpmatrix
\ \ \ 
\pmatrix{-10}&{10}&{5}&{0}\cr{10}&{-10}&{0}&{10}\cr
{5}&{0}&{-2}&{0}\cr{0}&{10}&{0}&{-2}\cr\endpmatrix$ . 

\noindent
$n=6:  \hskip20pt 
\pmatrix{6}&{0}&{-1}\cr{0}&{0}&{1}\cr{-1}&{1}&{0}\cr{2}&{2}&{-1}\cr
\endpmatrix 
\ \ \ 
\pmatrix{-12}&{12}&{6}&{0}\cr{12}&{-12}&{0}&{12}\cr
{6}&{0}&{-2}&{0}\cr{0}&{12}&{0}&{-4}\cr\endpmatrix$ . 

\noindent
$n=7: \hskip20pt  
\pmatrix{7}&{0}&{-1}\cr{0}&{0}&{1}\cr{-1}&{1}&{0}\cr{3}&{2}&{-1}\cr
\endpmatrix
\ \ \ 
\pmatrix{-14}&{14}&{7}&{0}\cr{14}&{-14}&{0}&{14}\cr{7}&{0}&{-2}&{1}\cr
{0}&{14}&{1}&{-2}\cr\endpmatrix$. 

\noindent
$n=8:\hskip20pt 
\pmatrix{8}&{0}&{-1}\cr{0}&{0}&{1}\cr{-1}&{1}&{0}\cr{8}&{8}&{-3}\cr
\endpmatrix
\ \ \ 
\pmatrix{-16}&{16}&{8}&{16}\cr{16}&{-16}&{0}&{48}\cr
{8}&{0}&{-2}&{0}\cr{16}&{48}&{0}&{-16}\cr\endpmatrix$. 

\noindent
$n=9: \hskip20pt 
\pmatrix{9}&{0}&{-1}\cr{0}&{0}&{1}\cr
{-1}&{1}&{0}\cr{4}&{2}&{-1}\cr\endpmatrix
\ \ \ 
\pmatrix{-18}&{18}&{9}&{0}\cr{18}&{-18}&{0}&{18}\cr
{9}&{0}&{-2}&{2}\cr{0}&{18}&{2}&{-2}\cr\endpmatrix$. 

\noindent
$n=10:\hskip20pt
\pmatrix{10}&{0}&{-1}\cr{0}&{0}&{1}\cr{-1}&{1}&{0}\cr{4}&{2}&{-1}\cr
\endpmatrix
\ \ \ 
\pmatrix{-20}&{20}&{10}&{0}\cr{20}&{-20}&{0}&{20}\cr{10}&{0}&{-2}&{2}\cr
{0}&{20}&{2}&{-4}\cr\endpmatrix$. 

\noindent
$n=11:\hskip20pt 
\pmatrix{11}&{0}&{-1}\cr{0}&{0}&{1}\cr{-1}&{1}&{0}\cr{33}&{33}&{-10}\cr
{88}&{66}&{-23}\cr{5}&{2}&{-1}\cr\endpmatrix
\ \ \ 
\pmatrix{-22}&{22}&{11}&{143}&{220}&{0}\cr{22}&{-22}&{0}&{220}&{506}&{22}\cr
{11}&{0}&{-2}&{0}&{22}&{3}\cr{143}&{220}&{0}&{-22}&{22}&{11}\cr
{220}&{506}&{22}&{22}&{-22}&{0}\cr{0}&{22}&{3}&{11}&{0}&{-2}\cr\endpmatrix$. 

\noindent
$n=12: \hskip20pt 
\pmatrix{12}&{0}&{-1}\cr{0}&{0}&{1}\cr{-1}&{1}&{0}\cr{3}&{3}&{-1}\cr
\endpmatrix
\ \ \ 
\pmatrix{-24}&{24}&{12}&{12}\cr{24}&{-24}&{0}&{24}\cr{12}&{0}&{-2}&{0}\cr
{12}&{24}&{0}&{-6}\cr\endpmatrix$. 

\noindent
$n=13: \hskip20pt 
\pmatrix{13}&{0}&{-1}\cr{0}&{0}&{1}\cr{-1}&{1}&{0}\cr{4}&{3}&{-1}\cr
{6}&{2}&{-1}\cr\endpmatrix
\ \ \ 
\pmatrix{-26}&{26}&{13}&{13}&{0}\cr
{26}&{-26}&{0}&{26}&{26}\cr{13}&{0}&{-2}&{1}&{4}\cr
{13}&{26}&{1}&{-2}&{0}\cr{0}&{26}&{4}&{0}&{-2}\cr\endpmatrix$ .

\noindent 
$n=14: \hskip20pt 
\pmatrix{14}&{0}&{-1}\cr{0}&{0}&{1}\cr{-1}&{1}&{0}\cr{7}&{7}&{-2}\cr
{6}&{2}&{-1}\cr\endpmatrix
\ \ \ 
\pmatrix{-28}&{28}&{14}&{42}&{0}\cr{28}&{-28}&{0}&{56}&{28}\cr
{14}&{0}&{-2}&{0}&{4}\cr{42}&{56}&{0}&{-14}&{0}\cr{0}&{28}&{4}&{0}&{-4}\cr
\endpmatrix$. 

\noindent 
$n=15: \hskip20pt 
\pmatrix{15}&{0}&{-1}\cr{0}&{0}&{1}\cr{-1}&{1}&{0}\cr
{15}&{15}&{-4}\cr{60}&{30}&{-11}\cr{7}&{2}&{-1}\cr\endpmatrix
\ \ \ 
\pmatrix{-30}&{30}&{15}&{105}&{120}&{0}\cr
{30}&{-30}&{0}&{120}&{330}&{30}\cr{15}&{0}&{-2}&{0}&{30}&{5}\cr
{105}&{120}&{0}&{-30}&{30}&{15}\cr{120}&{330}&{30}&{30}&{-30}&{0}\cr
{0}&{30}&{5}&{15}&{0}&{-2}\cr\endpmatrix$. 

\noindent
$n=16: \hskip20pt 
\pmatrix{16}&{0}&{-1}\cr{0}&{0}&{1}\cr{-1}&{1}&{0}\cr
{5}&{3}&{-1}\cr{48}&{16}&{-7}\cr\endpmatrix
\ \ \ 
\pmatrix{-32}&{32}&{16}&{16}&{32}\cr{32}&{-32}&{0}&{32}&{224}\cr
{16}&{0}&{-2}&{2}&{32}\cr{16}&{32}&{2}&{-2}&{0}\cr
{32}&{224}&{32}&{0}&{-32}\cr\endpmatrix$. 

\noindent 
$n=17: \hskip20pt 
\pmatrix{17}&{0}&{-1}\cr{0}&{0}&{1}\cr{-1}&{1}&{0}\cr{4}&{4}&{-1}\cr
{85}&{51}&{-16}\cr{204}&{102}&{-35}\cr{19}&{8}&{-3}\cr{8}&{2}&{-1}\cr
\endpmatrix 
\ \ \ 
\pmatrix{-34}&{34}&{17}&{34}&{323}&{544}&{34}&{0}\cr
{34}&{-34}&{0}&{34}&{544}&{1190}&{102}&{34}\cr
{17}&{0}&{-2}&{0}&{34}&{102}&{11}&{6}\cr
{34}&{34}&{0}&{-2}&{0}&{34}&{6}&{6}\cr
{323}&{544}&{34}&{0}&{-34}&{34}&{17}&{34}\cr
{544}&{1190}&{102}&{34}&{34}&{-34}&{0}&{34}\cr
{34}&{102}&{11}&{6}&{17}&{0}&{-2}&{0}\cr
{0}&{34}&{6}&{6}&{34}&{34}&{0}&{-2}\cr\endpmatrix 
$. 

\noindent 
$n=18: \hskip20pt 
\pmatrix{18}&{0}&{-1}\cr{0}&{0}&{1}\cr{-1}&{1}&{0}\cr
{4}&{4}&{-1}\cr{8}&{2}&{-1}\cr\endpmatrix
\ \ \ 
\pmatrix{-36}&{36}&{18}&{36}&{0}\cr{36}&{-36}&{0}&{36}&{36}\cr
{18}&{0}&{-2}&{0}&{6}\cr{36}&{36}&{0}&{-4}&{4}\cr
{0}&{36}&{6}&{4}&{-4}\cr\endpmatrix$. 

\noindent
$n=19: \hskip20pt 
\pmatrix{19}&{0}&{-1}\cr{0}&{0}&{1}\cr{-1}&{1}&{0}\cr{31}&{30}&{-7}\cr
{171}&{152}&{-37}\cr{190}&{152}&{-39}\cr{6}&{3}&{-1}\cr{9}&{2}&{-1}\cr
\endpmatrix
\ \ \ 
\pmatrix{-38}&{38}&{19}&{304}&{1482}&{1406}&{19}&{0}\cr
{38}&{-38}&{0}&{266}&{1406}&{1482}&{38}&{38}\cr
{19}&{0}&{-2}&{1}&{19}&{38}&{3}&{7}\cr
{304}&{266}&{1}&{-2}&{0}&{38}&{7}&{66}\cr
{1482}&{1406}&{19}&{0}&{-38}&{38}&{19}&{304}\cr
{1406}&{1482}&{38}&{38}&{38}&{-38}&{0}&{266}\cr
{19}&{38}&{3}&{7}&{19}&{0}&{-2}&{1}\cr
{0}&{38}&{7}&{66}&{304}&{266}&{1}&{-2}\cr\endpmatrix$. 

\noindent 
$n=20: \hskip20pt 
\pmatrix{20}&{0}&{-1}\cr{0}&{0}&{1}\cr
{-1}&{1}&{0}\cr{4}&{4}&{-1}\cr{15}&{5}&{-2}\cr\endpmatrix
\ \ \ 
\pmatrix{-40}&{40}&{20}&{40}&{20}\cr{40}&{-40}&{0}&{40}&{80}\cr
{20}&{0}&{-2}&{0}&{10}\cr{40}&{40}&{0}&{-8}&{0}\cr
{20}&{80}&{10}&{0}&{-10}\cr
\endpmatrix$  . 

\noindent 
$n=21: \hskip20pt 
\pmatrix{21}&{0}&{-1}\cr{0}&{0}&{1}\cr{-1}&{1}&{0}\cr{6}&{3}&{-1}\cr
{10}&{2}&{-1}\cr\endpmatrix 
\ \ \ 
\pmatrix{-42}&{42}&{21}&{21}&{0}\cr{42}&{-42}&{0}&{42}&{42}\cr
{21}&{0}&{-2}&{3}&{8}\cr{21}&{42}&{3}&{-6}&{0}\cr
{0}&{42}&{8}&{0}&{-2}\cr\endpmatrix$. 
 
\noindent 
$n=22: \hskip20pt 
\pmatrix{22}&{0}&{-1}\cr{0}&{0}&{1}\cr{-1}&{1}&{0}\cr
{14}&{14}&{-3}\cr{110}&{88}&{-21}\cr{132}&{88}&{-23}\cr
{7}&{3}&{-1}\cr{10}&{2}&{-1}\cr\endpmatrix 
\ \ \ 
\pmatrix{-44}&{44}&{22}&{176}&{1012}&{924}&{22}&{0}\cr
{44}&{-44}&{0}&{132}&{924}&{1012}&{44}&{44}\cr
{22}&{0}&{-2}&{0}&{22}&{44}&{4}&{8}\cr
{176}&{132}&{0}&{-4}&{0}&{44}&{8}&{36}\cr
{1012}&{924}&{22}&{0}&{-44}&{44}&{22}&{176}\cr
{924}&{1012}&{44}&{44}&{44}&{-44}&{0}&{132}\cr
{22}&{44}&{4}&{8}&{22}&{0}&{-2}&{0}\cr
{0}&{44}&{8}&{36}&{176}&{132}&{0}&{-4}\cr\endpmatrix$. 

\noindent 
$n=23$: Hyperbolic type with $(w,w)=-138$, 
the symmetry group generated by the translation $C$, and 
\newline 
$e=\pmatrix{23}&{0}&{-1}\cr{0}&{0}&{1}\cr
{-1}&{1}&{0}\cr{115}&{115}&{-24}\cr\endpmatrix\ \ \ 
G(e)=
\pmatrix{-46}&{46}&{23}&{1541}\cr{46}&{-46}&{0}&{1104}\cr
{23}&{0}&{-2}&{0}\cr{1541}&{1104}&{0}&{-46}\cr\endpmatrix$
\newline
$f=\pmatrix{368}&{138}&{-47}\cr{161}&{69}&{-22}\cr
{13}&{7}&{-2}\cr{276}&{184}&{-47}\cr\endpmatrix\ \ \ 
G(f)=
\pmatrix{-46}&{46}&{46}&{4186}\cr{46}&{-46}&{0}&{1104}\cr
{46}&{0}&{-2}&{0}\cr{4186}&{1104}&{0}&{-46}\cr\endpmatrix$
\newline 
$C=
\pmatrix{25}&{92}&{460}\cr{23}&{81}&{414}\cr{-5}&{-18}&{-91}\cr
\endpmatrix\ \  
w=\pmatrix{{-46}}\cr{{-23}}\cr{7}\cr\endpmatrix$. 

\noindent 
$n=24: \hskip20pt 
\pmatrix{24}&{0}&{-1}\cr{0}&{0}&{1}\cr
{-1}&{1}&{0}\cr{24}&{24}&{-5}\cr
{48}&{24}&{-7}\cr{144}&{48}&{-17}\cr
{19}&{5}&{-2}\cr{120}&{24}&{-11}\cr\endpmatrix
\ \ \ 
\pmatrix{-48}&{48}&{24}&{336}&{240}&{336}&{24}&{48}\cr
{48}&{-48}&{0}&{240}&{336}&{816}&{96}&{528}\cr
{24}&{0}&{-2}&{0}&{24}&{96}&{14}&{96}\cr
{336}&{240}&{0}&{-48}&{48}&{528}&{96}&{816}\cr
{240}&{336}&{24}&{48}&{-48}&{48}&{24}&{336}\cr
{336}&{816}&{96}&{528}&{48}&{-48}&{0}&{240}\cr
{24}&{96}&{14}&{96}&{24}&{0}&{-2}&{0}\cr
{48}&{528}&{96}&{816}&{336}&{240}&{0}&{-48}\cr\endpmatrix 
$. 

\noindent 
$n=25: \hskip20pt 
\pmatrix{25}&{0}&{-1}\cr
{0}&{0}&{1}\cr{-1}&{1}&{0}\cr
{6}&{4}&{-1}\cr{8}&{3}&{-1}\cr{12}&{2}&{-1}\cr\endpmatrix 
\ \ \ 
\pmatrix{-50}&{50}&{25}&{50}&{25}&{0}\cr
{50}&{-50}&{0}&{50}&{50}&{50}\cr
{25}&{0}&{-2}&{2}&{5}&{10}\cr
{50}&{50}&{2}&{-2}&{0}&{10}\cr{25}&{50}&{5}&{0}&{-2}&{2}\cr
{0}&{50}&{10}&{10}&{2}&{-2}\cr\endpmatrix $. 

\noindent 
$n=26: \hskip20pt 
\pmatrix{26}&{0}&{-1}\cr{0}&{0}&{1}\cr
{-1}&{1}&{0}\cr{6}&{4}&{-1}\cr
{208}&{78}&{-25}\cr{468}&{156}&{-53}\cr
{55}&{17}&{-6}\cr{12}&{2}&{-1}\cr\endpmatrix 
\ \ \ 
\pmatrix{-52}&{52}&{26}&{52}&{728}&{1300}&{130}&{0}\cr
{52}&{-52}&{0}&{52}&{1300}&{2756}&{312}&{52}\cr
{26}&{0}&{-2}&{2}&{130}&{312}&{38}&{10}\cr
{52}&{52}&{2}&{-4}&{0}&{52}&{10}&{8}\cr
{728}&{1300}&{130}&{0}&{-52}&{52}&{26}&{52}\cr
{1300}&{2756}&{312}&{52}&{52}&{-52}&{0}&{52}\cr
{130}&{312}&{38}&{10}&{26}&{0}&{-2}&{2}\cr
{0}&{52}&{10}&{8}&{52}&{52}&{2}&{-4}\cr\endpmatrix$.   

\noindent
$n=27$: Not reflective. 

\noindent 
$n=28: \hskip20pt  
\pmatrix{28}&{0}&{-1}\cr{0}&{0}&{1}\cr{-1}&{1}&{0}
\cr{21}&{21}&{-4}\cr{84}&{56}&{-13}\cr
{112}&{56}&{-15}\cr{9}&{3}&{-1}\cr{35}&{7}&{-3}\cr\endpmatrix
\ \ \ 
\pmatrix{-56}&{56}&{28}&{364}&{840}&{728}&{28}&{28}\cr
{56}&{-56}&{0}&{224}&{728}&{840}&{56}&{168}\cr
{28}&{0}&{-2}&{0}&{28}&{56}&{6}&{28}\cr
{364}&{224}&{0}&{-14}&{28}&{168}&{28}&{210}\cr
{840}&{728}&{28}&{28}&{-56}&{56}&{28}&{364}\cr
{728}&{840}&{56}&{168}&{56}&{-56}&{0}&{224}\cr
{28}&{56}&{6}&{28}&{28}&{0}&{-2}&{0}\cr
{28}&{168}&{28}&{210}&{364}&{224}&{0}&{-14}\cr\endpmatrix$.  

\noindent 
$n=29$: Hyperbolic type with $(w,w)=-348$, the symmetry group 
generated by the translation $C$, and  
\newline 
$e=\pmatrix{29}&{0}&{-1}\cr{0}&{0}&{1}\cr{-1}&{1}&{0}\cr{70}&{70}&{-13}\cr
{957}&{928}&{-175}\cr\endpmatrix \ \ \ 
G(e)=\pmatrix{-58}&{58}&{29}&{1276}&{16762}\cr
{58}&{-58}&{0}&{754}&{10150}\cr{29}&{0}&{-2}&{0}&{29}\cr
{1276}&{754}&{0}&{-2}&{0}\cr{16762}&{10150}&{29}&{0}&{-58}\cr\endpmatrix$
\newline   
$f=\pmatrix{580}&{174}&{-59}\cr{261}&{87}&{-28}\cr{26}&{10}&{-3}\cr
{7}&{4}&{-1}\cr{348}&{290}&{-59}\cr\endpmatrix \ \ \ 
G(f)=\pmatrix{-58}&{58}&{58}&{116}&{26854}\cr
{58}&{-58}&{0}&{29}&{10150}\cr{58}&{0}&{-2}&{0}&{754}\cr
{116}&{29}&{0}&{-2}&{0}\cr{26854}&{10150}&{754}&{0}&{-58}\cr\endpmatrix$ 
\newline 
$C=\pmatrix
{121}&{464}&{2552}\cr
{116}&{441}&{2436}\cr
{-22}&{-84}&{-463}\cr\endpmatrix\ \ \ 
w=\pmatrix{{-58}}\cr{{-29}}\cr{8}\cr\endpmatrix $. 

\noindent
$n=30: \hskip20pt 
\pmatrix{30}&{0}&{-1}\cr{0}&{0}&{1}\cr{-1}&{1}&{0}\cr
{5}&{5}&{-1}\cr{9}&{3}&{-1}\cr{14}&{2}&{-1}\cr\endpmatrix
\ \ \ 
\pmatrix{-60}&{60}&{30}&{90}&{30}&{0}\cr{60}&{-60}&{0}&{60}&{60}&{60}\cr
{30}&{0}&{-2}&{0}&{6}&{12}\cr{90}&{60}&{0}&{-10}&{0}&{20}\cr
{30}&{60}&{6}&{0}&{-6}&{0}\cr{0}&{60}&{12}&{20}&{0}&{-4}\cr\endpmatrix$ . 

\noindent 
$n=31$: Hyperbolic type with $(w,w)=-434$, the symmetry group 
generated by the translation $C$, and     
\newline  
$e=\pmatrix
{31}&{0}&{-1}\cr{0}&{0}&{1}\cr{-1}&{1}&{0}\cr
{6}&{5}&{-1}\cr{2728}&{1581}&{-373}\cr\endpmatrix \ \ \ 
G(e)=\pmatrix{-62}&{62}&{31}&{93}&{25885}\cr
{62}&{-62}&{0}&{62}&{23126}\cr
{31}&{0}&{-2}&{1}&{1147}\cr
{93}&{62}&{1}&{-2}&{0}\cr
{25885}&{23126}&{1147}&{0}&{-62}\cr\endpmatrix$
\newline  
$f=\pmatrix{3100}&{620}&{-249}\cr{744}&{155}&{-61}\cr
{10}&{3}&{-1}\cr{33}&{15}&{-4}\cr{496}&{248}&{-63}\cr\endpmatrix\ \ \ 
G(f)=
\pmatrix{-62}&{62}&{62}&{5208}&{103726}\cr{62}&{-62}&{0}&{1147}&{23126}\cr
{62}&{0}&{-2}&{1}&{62}\cr{5208}&{1147}&{1}&{-2}&{0}\cr
{103726}&{23126}&{62}&{0}&{-62}\cr\endpmatrix$
\newline  
$C=\pmatrix{484}&{2511}&{12276}\cr
{279}&{1444}&{7068}\cr{-66}&{-342}&{-1673}\cr\endpmatrix \ \ \ 
w=\pmatrix{{-93}}\cr{{-31}}\cr{10}\cr\endpmatrix$. 

\noindent
$n=32$: Not reflective.   

\noindent 
$n=33$: 
\newline  
$
\pmatrix{33}&{0}&{-1}\cr{0}&{0}&{1}\cr{-1}&{1}&{0}\cr
{11}&{11}&{-2}\cr{8}&{4}&{-1}\cr{99}&{33}&{-10}\cr
{264}&{66}&{-23}\cr{37}&{8}&{-3}\cr{121}&{22}&{-9}\cr{16}&{2}&{-1}\cr
\endpmatrix 
\ \ \ 
\pmatrix{-66}&{66}&{33}&{231}&{66}&{429}&{660}&{66}&{132}&{0}\cr
{66}&{-66}&{0}&{132}&{66}&{660}&{1518}&{198}&{594}&{66}\cr
{33}&{0}&{-2}&{0}&{4}&{66}&{198}&{29}&{99}&{14}\cr
{231}&{132}&{0}&{-22}&{0}&{132}&{594}&{99}&{385}&{66}\cr
{66}&{66}&{4}&{0}&{-2}&{0}&{66}&{14}&{66}&{14}\cr
{429}&{660}&{66}&{132}&{0}&{-66}&{66}&{33}&{231}&{66}\cr
{660}&{1518}&{198}&{594}&{66}&{66}&{-66}&{0}&{132}&{66}\cr
{66}&{198}&{29}&{99}&{14}&{33}&{0}&{-2}&{0}&{4}\cr
{132}&{594}&{99}&{385}&{66}&{231}&{132}&{0}&{-22}&{0}\cr
{0}&{66}&{14}&{66}&{14}&{66}&{66}&{4}&{0}&{-2}\cr\endpmatrix$. 

\noindent 
$n=34: \hskip20pt  
\pmatrix{34}&{0}&{-1}\cr{0}&{0}&{1}\cr{-1}&{1}&{0}\cr
{17}&{17}&{-3}\cr{8}&{4}&{-1}\cr{11}&{3}&{-1}\cr{16}&{2}&{-1}\cr
\endpmatrix 
\ \ \ 
\pmatrix{-68}&{68}&{34}&{374}&{68}&{34}&{0}\cr
{68}&{-68}&{0}&{204}&{68}&{68}&{68}\cr
{34}&{0}&{-2}&{0}&{4}&{8}&{14}\cr{374}&{204}&{0}&{-34}&{0}&{34}&{102}\cr
{68}&{68}&{4}&{0}&{-4}&{0}&{12}\cr{34}&{68}&{8}&{34}&{0}&{-2}&{2}\cr
{0}&{68}&{14}&{102}&{12}&{2}&{-4}\cr\endpmatrix$. 

\noindent 
$n=35$: Hyperbolic type with $(w,w)=-30$, the symmetry group 
generated by the central symmetry $C_1$ and the translation $C_2$, and      
\newline
$e=\pmatrix{35}&{0}&{-1}\cr{0}&{0}&{1}\cr
{-1}&{1}&{0}\cr{35}&{35}&{-6}\cr\endpmatrix\ \ \  
G(e)=\pmatrix
{-70}&{70}&{35}&{805}\cr{70}&{-70}&{0}&{420}\cr
{35}&{0}&{-2}&{0}\cr{805}&{420}&{0}&{-70}\cr\endpmatrix$
\newline 
$C_1=\pmatrix
{35}&{144}&{840}\cr{9}&{35}&{210}\cr{-3}&{-12}&{-71}\cr\endpmatrix\ \ \ 
C_2=\pmatrix
{7}&{45}&{210}\cr{5}&{28}&{140}\cr{-1}&{-6}&{-29}\cr\endpmatrix\ \ \ 
w=\pmatrix{{-30}}\cr{{-10}}\cr{3}\cr\endpmatrix$. 

\noindent 
$n=36: \hskip20pt 
\pmatrix{36}&{0}&{-1}\cr{0}&{0}&{1}\cr{-1}&{1}&{0}\cr{8}&{4}&{-1}\cr
{40}&{8}&{-3}\cr{63}&{9}&{-4}\cr\endpmatrix
\ \ \ 
\pmatrix{-72}&{72}&{36}&{72}&{72}&{36}\cr{72}&{-72}&{0}&{72}&{216}&{288}\cr
{36}&{0}&{-2}&{4}&{32}&{54}\cr{72}&{72}&{4}&{-8}&{8}&{36}\cr
{72}&{216}&{32}&{8}&{-8}&{0}\cr{36}&{288}&{54}&{36}&{0}&{-18}\cr
\endpmatrix$ .
 
\noindent
$n=37$: Hyperbolic type with $(w,w)=-66$, the symmetry group 
generated by the translation $C$, and 
\newline
$e=\pmatrix{37}&{0}&{-1}\cr{0}&{0}&{1}\cr
{-1}&{1}&{0}\cr{6}&{6}&{-1}\cr{259}&{185}&{-36}\cr\endpmatrix\ \ \ 
G(e)=\pmatrix{-74}&{74}&{37}&{148}&{4181}\cr
{74}&{-74}&{0}&{74}&{2664}\cr{37}&{0}&{-2}&{0}&{74}\cr
{148}&{74}&{0}&{-2}&{0}\cr{4181}&{2664}&{74}&{0}&{-74}\cr\endpmatrix$ 
\newline 
$f=\pmatrix{231}&{41}&{-16}\cr{12}&{3}&{-1}\cr\endpmatrix \ \ \ 
G(f)=\pmatrix{-2}&{1}\cr{1}&{-2}\cr\endpmatrix$
\newline 
$C=\pmatrix{37}&{225}&{1110}\cr
{25}&{148}&{740}\cr{-5}&{-30}&{-149}\cr\endpmatrix\ \ \  
w=\pmatrix{{-30}}\cr{{-10}}\cr{3}\cr\endpmatrix$ .

\noindent 
$n=38$: Hyperbolic type with $(w,w)=-532$, the symmetry 
group generated by the central symmetry $C_1$ and the 
translation $C_2$, and  
\newline
$e=\pmatrix{38}&{0}&{-1}\cr{0}&{0}&{1}\cr{-1}&{1}&{0}\cr
{6}&{6}&{-1}\cr{304}&{190}&{-39}\cr\endpmatrix\ \ \ 
G(e)=\pmatrix{-76}&{76}&{38}&{152}&{4256}\cr
{76}&{-76}&{0}&{76}&{2964}\cr{38}&{0}&{-2}&{0}&{114}\cr
{152}&{76}&{0}&{-4}&{0}\cr{4256}&{2964}&{114}&{0}&{-76}\cr\endpmatrix$
\newline 
$C_1=\pmatrix{38}&{169}&{988}\cr
{9}&{38}&{228}\cr{-3}&{-13}&{-77}\cr\endpmatrix 
\ \ \ 
C_2=\pmatrix
{32}&{171}&{912}\cr{19}&{98}&{532}\cr
{-4}&{-21}&{-113}\cr\endpmatrix\ \ \ 
w=\pmatrix{{-114}}\cr{{-38}}\cr{11}\cr\endpmatrix$.
  
\noindent 
$n=39: \hskip20pt 
\pmatrix{39}&{0}&{-1}\cr{0}&{0}&{1}\cr{-1}&{1}&{0}\cr
{6}&{6}&{-1}\cr{26}&{13}&{-3}\cr{12}&{3}&{-1}\cr{19}&{2}&{-1}\cr
\endpmatrix 
\ \ \ 
\pmatrix{-78}&{78}&{39}&{156}&{273}&{39}&{0}\cr
{78}&{-78}&{0}&{78}&{234}&{78}&{78}\cr
{39}&{0}&{-2}&{0}&{13}&{9}&{17}\cr{156}&{78}&{0}&{-6}&{0}&{12}&{48}\cr
{273}&{234}&{13}&{0}&{-26}&{0}&{65}\cr{39}&{78}&{9}&{12}&{0}&{-6}&{3}\cr
{0}&{78}&{17}&{48}&{65}&{3}&{-2}\cr\endpmatrix$ . 

\noindent
$n=40$: Hyperbolic type with $(w,w)=-120$, the symmetry group  
generated by the skew-symmetry $C$, and    
\newline  
$e=\pmatrix{40}&{0}&{-1}\cr{0}&{0}&{1}\cr{-1}&{1}&{0}\cr
{120}&{120}&{-19}\cr{280}&{240}&{-41}\cr\endpmatrix\ \ \ 
G(e)=\pmatrix
{-80}&{80}&{40}&{3280}&{6320}\cr
{80}&{-80}&{0}&{1520}&{3280}\cr{40}&{0}&{-2}&{0}&{40}\cr
{3280}&{1520}&{0}&{-80}&{80}\cr{6320}&{3280}&{40}&{80}&{-80}\cr\endpmatrix$
\newline 
$C=\pmatrix{32}&{45}&{480}\cr{5}&{8}&{80}\cr{-2}&{-3}&{-31}\cr\endpmatrix
\ \ \ 
w=\pmatrix{{-30}}\cr{{-10}}\cr{3}\cr\endpmatrix$. 

\noindent 
$n=41$: Not reflective.

\noindent 
$n=42$:
\newline   
$\pmatrix{42}&{0}&{-1}\cr{0}&{0}&{1}\cr{-1}&{1}&{0}\cr{6}&{6}&{-1}\cr
{10}&{4}&{-1}\cr{168}&{42}&{-13}\cr{420}&{84}&{-29}\cr
{61}&{11}&{-4}\cr{114}&{18}&{-7}\cr{20}&{2}&{-1}\cr\endpmatrix
\pmatrix{-84}&{84}&{42}&{168}&{84}&{672}&{1092}&{126}&{168}&{0}\cr
{84}&{-84}&{0}&{84}&{84}&{1092}&{2436}&{336}&{588}&{84}\cr
{42}&{0}&{-2}&{0}&{6}&{126}&{336}&{50}&{96}&{18}\cr
{168}&{84}&{0}&{-12}&{0}&{168}&{588}&{96}&{204}&{48}\cr
{84}&{84}&{6}&{0}&{-4}&{0}&{84}&{18}&{48}&{16}\cr
{672}&{1092}&{126}&{168}&{0}&{-84}&{84}&{42}&{168}&{84}\cr
{1092}&{2436}&{336}&{588}&{84}&{84}&{-84}&{0}&{84}&{84}\cr
{126}&{336}&{50}&{96}&{18}&{42}&{0}&{-2}&{0}&{6}\cr
{168}&{588}&{96}&{204}&{48}&{168}&{84}&{0}&{-12}&{0}\cr
{0}&{84}&{18}&{48}&{16}&{84}&{84}&{6}&{0}&{-4}\cr\endpmatrix$. 

\noindent 
$n=43$: Hyperbolic type with $(w,w)=-94$, the symmetry group 
generated by the skew-symmetry $C$, and 
\newline
$e=\pmatrix
{21}&{2}&{-1}\cr{43}&{0}&{-1}\cr{0}&{0}&{1}\cr
{-1}&{1}&{0}\cr{7}&{6}&{-1}\cr\endpmatrix\ \ \ 
G(e)=\pmatrix{-2}&{0}&{86}&{19}&{54}\cr{0}&{-86}&{86}&{43}&{172}\cr
{86}&{86}&{-86}&{0}&{86}\cr{19}&{43}&{0}&{-2}&{1}\cr
{54}&{172}&{86}&{1}&{-2}\cr\endpmatrix $
\newline  
$C=\pmatrix{172}&{289}&{2924}\cr
{25}&{43}&{430}\cr{-10}&{-17}&{-171}\cr\endpmatrix\ \ \ 
w=\pmatrix{{-34}}\cr{{-10}}\cr{3}\cr\endpmatrix$. 

\noindent 
$n=44$: Hyperbolic type with $(w,w)=-66$, the symmetry group 
generated by the central symmetry $C_1$ and the  
translation $C_2$, and 
\newline 
$e=\pmatrix{44}&{0}&{-1}\cr{0}&{0}&{1}\cr{-1}&{1}&{0}\cr
{33}&{33}&{-5}\cr{176}&{132}&{-23}\cr\endpmatrix\ \ \ 
G(e)=\pmatrix{-88}&{88}&{44}&{1012}&{3784}\cr
{88}&{-88}&{0}&{440}&{2024}\cr{44}&{0}&{-2}&{0}&{44}\cr
{1012}&{440}&{0}&{-22}&{44}\cr{3784}&{2024}&{44}&{44}&{-88}\cr\endpmatrix$
\newline 
$C_1=\pmatrix{539}&{2916}&{16632}\cr{100}&{539}&{3080}\cr 
{-35}&{-189}&{-1079}\cr\endpmatrix \ \ \ 
C_2=\pmatrix
{16}&{99}&{528}\cr{11}&{64}&{352}\cr{-2}&{-12}&{-65}\cr\endpmatrix\ \  
w=\pmatrix{{-33}}\cr{{-11}}\cr{3}\cr\endpmatrix$.  
 
\noindent 
$n=45: \hskip10pt      
\left(\smallmatrix
{45}&{0}&{-1}\cr{0}&{0}&{1}\cr{-1}&{1}&{0}\cr{20}&{20}&{-3}\cr
{55}&{40}&{-7}\cr{58}&{38}&{-7}\cr{225}&{135}&{-26}\cr
{180}&{90}&{-19}\cr{11}&{4}&{-1}\cr{40}&{10}&{-3}\cr
{35}&{5}&{-2}\cr{22}&{2}&{-1}\cr\endsmallmatrix\right) \ \ \  
\left(\smallmatrix 
{-90}&{90}&{45}&{630}&{1170}&{1080}&{3735}&{2340}&{90}&{180}&{45}&{0}\cr
{90}&{-90}&{0}&{270}&{630}&{630}&{2340}&{1710}&{90}&{270}&{180}&{90}\cr
{45}&{0}&{-2}&{0}&{15}&{20}&{90}&{90}&{7}&{30}&{30}&{20}\cr
{630}&{270}&{0}&{-10}&{10}&{30}&{180}&{270}&{30}&{190}&{260}&{210}\cr
{1170}&{630}&{15}&{10}&{-10}&{0}&{45}&{180}&{30}&{260}&{415}&{360}\cr
{1080}&{630}&{20}&{30}&{0}&{-2}&{0}&{90}&{20}&{210}&{360}&{322}\cr
{3735}&{2340}&{90}&{180}&{45}&{0}&{-90}&{90}&{45}&{630}&{1170}&{1080}\cr
{2340}&{1710}&{90}&{270}&{180}&{90}&{90}&{-90}&{0}&{270}&{630}&{630}\cr
{90}&{90}&{7}&{30}&{30}&{20}&{45}&{0}&{-2}&{0}&{15}&{20}\cr
{180}&{270}&{30}&{190}&{260}&{210}&{630}&{270}&{0}&{-10}&{10}&{30}\cr
{45}&{180}&{30}&{260}&{415}&{360}&{1170}&{630}&{15}&{10}&{-10}&{0}\cr
{0}&{90}&{20}&{210}&{360}&{322}&{1080}&{630}&{20}&{30}&{0}&{-2}\cr
\endsmallmatrix\right)$. 

\noindent 
$n=46$: Hyperbolic type with $(w,w)=-184$, the symmetry group 
generated by the central symmetry $C_1$ and by the 
translation $C_2$, and 
\newline    
$e=\pmatrix{46}&{0}&{-1}\cr{0}&{0}&{1}\cr{-1}&{1}&{0}\cr
{529}&{529}&{-78}\cr{158}&{154}&{-23}\cr{966}&{920}&{-139}\cr\endpmatrix
\ \ \ 
G(e)=\pmatrix{-92}&{92}&{46}&{17158}&{4968}&{29532}\cr
{92}&{-92}&{0}&{7176}&{2116}&{12788}\cr{46}&{0}&{-2}&{0}&{4}&{46}\cr
{17158}&{7176}&{0}&{-46}&{0}&{230}\cr{4968}&{2116}&{4}&{0}&{-4}&{0}\cr
{29532}&{12788}&{46}&{230}&{0}&{-92}\cr\endpmatrix$
\newline   
$C_1=\pmatrix
{23}&{32}&{368}\cr{18}&{23}&{276}\cr{-3}&{-4}&{-47}\cr\endpmatrix\ \ \ 
C_2= 
\pmatrix
{49}&{184}&{1288}\cr{46}&{169}&{1196}\cr{-7}&{-26}&{-183}\cr\endpmatrix\ \ \ 
w=\pmatrix
{{-46}}\cr{{-23}}\cr{5}\cr\endpmatrix$ . 
 
\noindent
$n=47$: Not reflective. 

\noindent 
$n=48$: Hyperbolic type with $(w,w)=-352$, the symmetry 
group generated by the translation $C$, and  
\newline 
$e=\pmatrix{48}&{0}&{-1}\cr{0}&{0}&{1}\cr{-1}&{1}&{0}\cr
{48}&{48}&{-7}\cr{144}&{96}&{-17}\cr\endpmatrix \ \ \ 
G(e)=\pmatrix{-96}&{96}&{48}&{1632}&{2976}\cr
{96}&{-96}&{0}&{672}&{1632}\cr{48}&{0}&{-2}&{0}&{48}\cr
{1632}&{672}&{0}&{-96}&{96}\cr{2976}&{1632}&{48}&{96}&{-96}\cr\endpmatrix$,
\newline 
$f=\pmatrix{285}&{33}&{-14}\cr{15}&{3}&{-1}\cr\endpmatrix\ \ \  
G(f)=\pmatrix{-6}&{6}\cr{6}&{-6}\cr\endpmatrix$,
\newline 
$C=\pmatrix{27}&{256}&{1152}\cr
{16}&{147}&{672}\cr{-3}&{-28}&{-127}\cr\endpmatrix\ \ \ 
w=\pmatrix{{-64}}\cr{{-16}}\cr{5}\cr\endpmatrix$. 

\noindent 
$n=49$:
\newline   
$\left(\smallmatrix{49}&{0}&{-1}\cr{0}&{0}&{1}\cr{-1}&{1}&{0}\cr
{8}&{6}&{-1}\cr{104}&{57}&{-11}\cr{451}&{240}&{-47}\cr
{396}&{208}&{-41}\cr{1421}&{735}&{-146}\cr{980}&{490}&{-99}\cr
{51}&{24}&{-5}\cr{12}&{4}&{-1}\cr{16}&{3}&{-1}\cr{39}&{5}&{-2}\cr
{24}&{2}&{-1}\cr\endsmallmatrix\right) \ \ 
\left(\smallmatrix
{-98}&{98}&{49}&{196}&{1715}&{7154}&{6174}&{21707}&{14308}&
{686}&{98}&{49}&{49}&{0}\cr
{98}&{-98}&{0}&{98}&{1078}&{4606}&{4018}&{14308}&{9702}&
{490}&{98}&{98}&{196}&{98}\cr
{49}&{0}&{-2}&{2}&{47}&{211}&{188}&{686}&{490}&{27}&{8}&{13}&{34}&{22}\cr
{196}&{98}&{2}&{-2}&{2}&{20}&{22}&{98}&{98}&{8}&{6}&{22}&{78}&{62}\cr
{1715}&{1078}&{47}&{2}&{-2}&{1}&{6}&{49}&{98}&{13}&{22}&{146}&{587}&{498}\cr
{7154}&{4606}&{211}&{20}&{1}&{-2}&{2}&{49}&{196}&{34}&{78}&{587}&{2403}
&{2056}\cr
{6174}&{4018}&{188}&{22}&{6}&{2}&{-2}&{0}&{98}&{22}&{62}&{498}&{2056}&{1766}\cr
{21707}&{14308}&{686}&{98}&{49}&{49}&{0}&{-98}&{98}&{49}&{196}&{1715}
&{7154}&{6174}\cr
{14308}&{9702}&{490}&{98}&{98}&{196}&{98}&{98}&{-98}&{0}&{98}&{1078}&
{4606}&{4018}\cr
{686}&{490}&{27}&{8}&{13}&{34}&{22}&{49}&{0}&{-2}&{2}&{47}&{211}&{188}\cr
{98}&{98}&{8}&{6}&{22}&{78}&{62}&{196}&{98}&{2}&{-2}&{2}&{20}&{22}\cr
{49}&{98}&{13}&{22}&{146}&{587}&{498}&{1715}&{1078}&{47}&{2}&{-2}&{1}&{6}\cr
{49}&{196}&{34}&{78}&{587}&{2403}&{2056}&{7154}&{4606}&{211}&
{20}&{1}&{-2}&{2}\cr{0}&{98}&{22}&{62}&{498}&{2056}&{1766}&{6174}&
{4018}&{188}&{22}&{6}&{2}&{-2}\cr\endsmallmatrix\right)$.

\noindent 
$n=50:  \left(\smallmatrix{50}&{0}&{-1}\cr{0}&{0}&{1}\cr
{-1}&{1}&{0}\cr{8}&{6}&{-1}\cr
{12}&{4}&{-1}\cr{78}&{16}&{-5}\cr{800}&{150}&{-49}\cr{1700}&{300}&{-101}\cr
{239}&{41}&{-14}\cr{434}&{72}&{-25}\cr
{56}&{8}&{-3}\cr{24}&{2}&{-1}\cr\endsmallmatrix\right)\ \ \  
\left(\smallmatrix
{-100}&{100}&{50}&{200}&{100}&{300}&{2600}&{4900}&{650}&{1100}&{100}&{0}\cr
{100}&{-100}&{0}&{100}&{100}&{500}&{4900}&{10100}&{1400}&{2500}&{300}&{100}\cr
{50}&{0}&{-2}&{2}&{8}&{62}&{650}&{1400}&{198}&{362}&{48}&{22}\cr
{200}&{100}&{2}&{-4}&{4}&{96}&{1100}&{2500}&{362}&{680}&{100}&{60}\cr
{100}&{100}&{8}&{4}&{-4}&{4}&{100}&{300}&{48}&{100}&{20}&{20}\cr
{300}&{500}&{62}&{96}&{4}&{-4}&{0}&{100}&{22}&{60}&{20}&{40}\cr
{2600}&{4900}&{650}&{1100}&{100}&{0}&{-100}&{100}&{50}&{200}&{100}&{300}\cr
{4900}&{10100}&{1400}&{2500}&{300}&{100}&{100}&{-100}&{0}&{100}&{100}&{500}\cr
{650}&{1400}&{198}&{362}&{48}&{22}&{50}&{0}&{-2}&{2}&{8}&{62}\cr
{1100}&{2500}&{362}&{680}&{100}&{60}&{200}&{100}&{2}&{-4}&{4}&{96}\cr
{100}&{300}&{48}&{100}&{20}&{20}&{100}&{100}&{8}&{4}&{-4}&{4}\cr
{0}&{100}&{22}&{60}&{20}&{40}&{300}&{500}&{62}&{96}&{4}&{-4}\cr
\endsmallmatrix\right)$. 

\noindent 
$n=51$: Not reflective.  
 
\noindent 
$n=52$: Hyperbolic type with $(w,w)=-1040$, the symmetry group 
generated by the skew-symmetry $C$, and 
\newline 
$e=\pmatrix{491}&{61}&{-24}\cr{108}&{12}&{-5}\cr
{143}&{13}&{-6}\cr{52}&{0}&{-1}\cr
{0}&{0}&{1}\cr{-1}&{1}&{0}\cr\endpmatrix\ \ \ 
G(e)=\pmatrix{-2}&{0}&{130}&{676}&{2496}&{430}\cr
{0}&{-8}&{0}&{104}&{520}&{96}\cr{130}&{0}&{-26}&{52}&{624}&{130}\cr
{676}&{104}&{52}&{-104}&{104}&{52}\cr
{2496}&{520}&{624}&{104}&{-104}&{0}\cr
{430}&{96}&{130}&{52}&{0}&{-2}\cr\endpmatrix$ 
\newline 
$C=\pmatrix{100}&{117}&{1560}\cr
{13}&{16}&{208}\cr{-5}&{-6}&{-79}\cr\endpmatrix\ \ \ 
w=\pmatrix{{-78}}\cr{{-26}}\cr{7}\cr\endpmatrix$ . 

\noindent
$n=53$: Not reflective. 

\noindent 
$n=54$: Not reflective.  

\noindent
$n=55$:
\newline  
$\left(\smallmatrix
{55}&{0}&{-1}\cr{0}&{0}&{1}\cr{-1}&{1}&{0}\cr{22}&{22}&{-3}\cr
{10}&{5}&{-1}\cr{44}&{11}&{-3}\cr{18}&{3}&{-1}\cr{440}&{55}&{-21}\cr
{1980}&{220}&{-89}\cr{361}&{39}&{-16}\cr{3938}&{418}&{-173}\cr
{485}&{50}&{-21}\cr{451}&{44}&{-19}\cr{27}&{2}&{-1}\cr
\endsmallmatrix\right) \ \ 
\left(\smallmatrix
{-110}&{110}&{55}&{880}&{165}&{275}&{55}&{715}&{2310}
&{385}&{3960}&{440}&{330}&{0}\cr
{110}&{-110}&{0}&{330}&{110}&{330}&{110}&{2310}&{9790}&
{1760}&{19030}&{2310}&{2090}&{110}\cr
{55}&{0}&{-2}&{0}&{5}&{33}&{15}&{385}&{1760}&{322}&{3520}&{435}&{407}&{25}\cr
{880}&{330}&{0}&{-22}&{0}&{220}&{132}&{3960}&{19030}&
{3520}&{38742}&{4840}&{4620}&{308}\cr
{165}&{110}&{5}&{0}&{-10}&{0}&{10}&{440}&{2310}&{435}&{4840}
&{615}&{605}&{45}\cr
{275}&{330}&{33}&{220}&{0}&{-22}&{0}&{330}&{2090}&{407}
&{4620}&{605}&{627}&{55}\cr
{55}&{110}&{15}&{132}&{10}&{0}&{-2}&{0}&{110}&{25}&{308}&{45}&{55}&{7}\cr
{715}&{2310}&{385}&{3960}&{440}&{330}&{0}&{-110}&
{110}&{55}&{880}&{165}&{275}&{55}\cr
{2310}&{9790}&{1760}&{19030}&{2310}&{2090}&{110}&{110}&{-110}&
{0}&{330}&{110}&{330}&{110}\cr
{385}&{1760}&{322}&{3520}&{435}&{407}&{25}&{55}&{0}&{-2}&{0}&{5}&{33}&{15}\cr
{3960}&{19030}&{3520}&{38742}&{4840}&{4620}&{308}&
{880}&{330}&{0}&{-22}&{0}&{220}&{132}\cr
{440}&{2310}&{435}&{4840}&{615}&{605}&{45}&{165}
&{110}&{5}&{0}&{-10}&{0}&{10}\cr
{330}&{2090}&{407}&{4620}&{605}&{627}&{55}&{275}&{330}&{33}&
{220}&{0}&{-22}&{0}\cr
{0}&{110}&{25}&{308}&{45}&{55}&{7}&{55}&{110}&{15}&{132}&{10}&{0}&{-2}\cr
\endsmallmatrix\right)$.

$n=56$: Hyperbolic type with $(w,w)=-56$, the symmetry group  
generated by the central symmetry $C_1$ and the 
translation $C_2$, and 
\newline  
$e=\pmatrix{56}&{0}&{-1}\cr{0}&{0}&{1}\cr{-1}&{1}&{0}\cr{7}&{7}&{-1}\cr
{11}&{5}&{-1}\cr{168}&{56}&{-13}\cr{224}&{56}&{-15}\cr\endpmatrix
\ \ \ 
G(e)=\pmatrix{-112}&{112}&{56}&{280}&{168}&{1680}&{1456}\cr
{112}&{-112}&{0}&{112}&{112}&{1456}&{1680}\cr
{56}&{0}&{-2}&{0}&{6}&{112}&{168}\cr
{280}&{112}&{0}&{-14}&{0}&{112}&{280}\cr
{168}&{112}&{6}&{0}&{-2}&{0}&{56}\cr
{1680}&{1456}&{112}&{112}&{0}&{-112}&{112}\cr
{1456}&{1680}&{168}&{280}&{56}&{112}&{-112}\cr\endpmatrix$
\newline
$C_1=\pmatrix{1183}&{10952}&{53872}\cr{128}&{1183}&{5824}\cr
{-52}&{-481}&{-2367}\cr\endpmatrix \ \ \ 
C_2=\pmatrix{32}&{343}&{1568}\cr{7}&{72}&{336}\cr
{-2}&{-21}&{-97}\cr\endpmatrix \ \ 
w=\pmatrix{{-98}}\cr{{-14}}\cr{5}\cr\endpmatrix$ . 

\noindent 
$n=57$: Hyperbolic type with $(w,w)=-494$, the symmetry group 
generated by the skew-symmetry $C$, and 
\newline 
$e=\pmatrix{321}&{30}&{-13}\cr{28}&{2}&{-1}\cr{57}&{0}&{-1}\cr
{0}&{0}&{1}\cr{-1}&{1}&{0}\cr{9}&{6}&{-1}\cr\endpmatrix\ \ \ 
G(e)=\pmatrix{-6}&{0}&{228}&{1482}&{291}&{714}\cr
{0}&{-2}&{0}&{114}&{26}&{72}\cr{228}&{0}&{-114}&{114}&{57}&{228}\cr
{1482}&{114}&{114}&{-114}&{0}&{114}\cr{291}&{26}&{57}&{0}&{-2}&{3}\cr
{714}&{72}&{228}&{114}&{3}&{-6}\cr\endpmatrix$
\newline 
$C=\pmatrix{48}&{475}&{2280}\cr
{19}&{192}&{912}\cr{-4}&{-40}&{-191}\cr\endpmatrix \ \ \ 
w=\pmatrix{{-95}}\cr{{-19}}\cr{6}\cr\endpmatrix$.

$n=58$: Not reflective.

$n=59$: Not reflective. 

$n=60$: Hyperbolic type with $(w,w)=-80$, the 
symmetry group generated by the central symmetry $C_1$ and 
the translation $C_2$, and 
\newline 
$e=\pmatrix{60}&{0}&{-1}\cr{0}&{0}&{1}\cr{-1}&{1}&{0}\cr
{15}&{15}&{-2}\cr{120}&{60}&{-11}\cr\endpmatrix\ \ \ 
G(e)=\pmatrix{-120}&{120}&{60}&{660}&{2280}\cr
{120}&{-120}&{0}&{240}&{1320}\cr{60}&{0}&{-2}&{0}&{60}\cr
{660}&{240}&{0}&{-30}&{60}\cr{2280}&{1320}&{60}&{60}&{-120}\cr\endpmatrix$
\newline 
$C_1=\pmatrix{12615}&{110224}&{577680}\cr{1444}&{12615}&{66120}\cr
{-551}&{-4814}&{-25231}\cr\endpmatrix\ \ \ 
C_2=\pmatrix{12}&{125}&{600}\cr
{5}&{48}&{240}\cr{-1}&{-10}&{-49}\cr\endpmatrix\ \ \ 
w=\pmatrix{{-50}}\cr{{-10}}\cr{3}\cr\endpmatrix$. 

As a result we get that 
cases 
$$
k=1 \text{\ ---\ } 22, 24 \text{\ ---\ }26, 
28,\,30,\,33,\,34,\,36,\,39,\, 42,\, 45,\,49,\,50,55
$$
are reflective of elliptic type; cases 
$$
k=23,\,29,\,31,\,35,\,37,\,38,\,40,\,43,\,44,\,46,\,48,\,52,\,56,\,60
$$
are reflective of hyperbolic type, and cases 
$$
k=27,\,32,\,41,\,47,\,51,\,53,\,54,\, 58,\, 59
$$
are not reflective. 
We can see that there are a lot of cases of hyperbolic type.

\smallpagebreak 
  
These calculations are mirror symmetric to results in \cite{G}, 
\cite{GH}, \cite{GN1}---\cite{GN6} where this series was considered for 
the model B.


\Refs 
\widestnumber\key{vedG2}


\ref
\key B1 
\by R. Borcherds
\paper Generalized Kac--Moody algebras
\jour J. of Algebra
\vol 115
\yr 1988
\pages 501--512
\endref

\ref
\key B2 
\by R. Borcherds
\paper The monster Lie algebra
\jour Adv. Math.
\vol 83
\yr 1990
\pages 30--47
\endref
\ref
\key B3
\by R. Borcherds
\paper The monstrous moonshine and monstrous Lie superalgebras
\jour Invent. Math.
\vol 109
\yr 1992
\pages 405--444
\endref


\ref
\key B4 
\by R. Borcherds
\paper Sporadic groups and string theory
\inbook Proc. European Congress of Mathem. 1992
\pages 411--421
\endref


\ref
\key B5 
\by R. Borcherds
\paper Automorphic forms on $O_{s+2,2}$ and
infinite products
\jour Invent. Math. \vol 120
\yr 1995
\pages 161--213
\endref

\ref
\key B6 
\by R. Borcherds
\paper The moduli space of Enriques surfaces and the fake monster Lie
superalgebra
\jour Topology 
\yr 1996 
\vol 35 \issue 3 
\pages 699--710 
\endref

\ref\key G 
\by V.A. Gritsenko
\paper Modular forms and moduli spaces of Abelian and K3 surfaces
\jour Algebra i Analyz
\vol 6:6
\yr 1994
\pages 65--102
\transl\nofrills  English transl. in
\jour St.Petersburg Math. Jour.
\vol 6:6
\yr 1995
\pages 1179--1208
\endref

\ref 
\key GH
\by V. Gritsenko, K. Hulek
\paper Commutator coverings of Siegel threefolds 
\jour Preprint RIMS Kyoto University 
\vol RIMS-1128 
\yr 1997
\moreref alg-geom/9702007
\endref

\ref
\key GN1 
\by V.A. Gritsenko, V.V. Nikulin
\paper Siegel automorphic form correction of some Lorentzi\-an
Kac--Moody Lie algebras
\jour Amer. J. Math.
\yr 1997 \toappear
\moreref alg-geom/9504006 
\endref

\ref
\key GN2 
\by V.A. Gritsenko, V.V. Nikulin
\paper Siegel automorphic form correction of a Lorentzian
Kac--Moody algebra
\jour C. R. Acad. Sci. Paris S\'er. A--B
\vol 321
\yr 1995
\pages 1151--1156
\endref

\ref
\key GN3 
\by V.A. Gritsenko, V.V. Nikulin
\paper K3 surfaces, Lorentzian Kac--Moody algebras and
mirror symmetry
\jour  Math. Res. Lett. \yr 1996 \vol 3 \issue 2 \pages 211--229 
\moreref  alg-geom/9510008
 \endref

\ref
\key GN4 
\by V.A. Gritsenko, V.V. Nikulin
\paper The Igusa modular forms and ``the simplest''
Lorentzian Kac--Moody algebras
\jour Matem. Sbornik 
\yr 1996 \vol 187 \issue 11 
\moreref alg-geom/9603010 
\endref

\ref
\key GN5 
\by V.A. Gritsenko, V.V. Nikulin
\paper Automorphic forms and Lorentzian Kac-Moody algebras.
Part I 
\jour Preprint RIMS Kyoto Univ. \yr 1996 
\vol RIMS-1116 
\moreref alg-geom/9610022
\endref

\ref
\key GN6 
\by V.A. Gritsenko, V.V. Nikulin
\paper Automorphic forms and Lorentzian Kac-Moody algebras.
Part II 
\jour Preprint RIMS Kyoto Univ. 
\yr 1996   
\vol RIMS-1122  
\moreref alg-geom/9611028 
\endref

\ref
\key GN7
\by V.A. Gritsenko, V.V. Nikulin
\paper The arithmetic mirror symmetry and Calabi--Yau manifolds 
\jour Preprint RIMS Kyoto Univ.    
\yr 1997 
\vol RIMS-1129 
\moreref alg-geom/9612002 
\endref

\ref
\key K1 
\by V. Kac
\book Infinite dimensional Lie algebras
\yr 1990
\publ Cambridge Univ. Press
\endref

\ref
\key Ku
\by Vic. S. Kulikov
\paper Degenerations of K3 surfaces and Enriques surfaces
\jour Izv. Akad. Nauk SSSR Ser. Mat.
\vol  41  \yr 1977 \pages 1008--1042
\transl\nofrills English transl. in
\jour Math. USSR Izv.
\vol 11 \yr 1977
\endref

\ref
\key N1 
\by V.V. Nikulin
\paper Finite automorphism groups of K\"ahler K3 surfaces
\jour Trudy Moskov. Mat. Obshch.
\vol 37
\yr 1979 \pages 73--137
\transl\nofrills English transl. in
\jour Trans. Moscow Math. Soc.
\vol 38 \issue 2 \yr 1980
\endref

\ref
\key N2 
\by V.V. Nikulin
\paper Integral symmetric bilinear forms and some of
their geometric applications
\jour Izv. Akad. Nauk SSSR Ser. Mat.
\vol  43
\yr 1979
\pages 111--177
\transl\nofrills English transl. in
\jour Math. USSR Izv.
\vol 14
\yr 1980
\endref

\ref
\key N3 
\by V.V. Nikulin
\paper On factor groups of the automorphism groups of
hyperbolic forms modulo subgroups generated by 2-reflections 
\jour Dokl. Akad. Nauk SSSR 
\yr 1979 
\vol 248 
\pages 1307--1309 
\transl\nofrills English transl. in
\jour Soviet Math. Dokl. 
\yr 1979 
\vol 20 
\pages 1156--1158 
\endref

\ref
\key N4 
\by V.V. Nikulin
\paper On the quotient groups of the automorphism groups of
hyperbolic forms by the subgroups generated by 2-reflections,
Algebraic-geometric applications
\jour Current Problems in Math. Vsesoyuz. Inst. Nauchn. i
Tekhn. Informatsii, Moscow
\yr 1981 \vol 18 
\pages 3--114
\transl\nofrills English transl. in
\jour J. Soviet Math.
\yr 1983
\vol 22
\pages 1401--1476
\endref

\ref
\key N5 
\by V.V. Nikulin
\paper On arithmetic groups generated by
reflections in Lobachevsky spaces
\jour Izv. Akad. Nauk SSSR Ser. Mat.
\vol  44   \yr 1980 \pages 637--669
\transl\nofrills English transl. in
\jour Math. USSR Izv.
\vol 16 \yr 1981
\endref

\ref
\key N6 
\by V.V. Nikulin
\paper On the classification of arithmetic groups generated by
reflections in Lobachevsky spaces
\jour Izv. Akad. Nauk SSSR Ser. Mat.
\vol  45
\issue 1
\yr 1981
\pages 113--142
\transl\nofrills English transl. in
\jour Math. USSR Izv.
\vol 18
\yr 1982
\endref

\ref
\key N7 
\by V.V. Nikulin
\paper Involutions of integral quadratic forms and their
applications to real algebraic geometry
\jour Izv. Akad. Nauk SSSR Ser. Mat.
\vol  47 \issue 1  \yr 1983
\transl\nofrills English transl. in
\jour Math. USSR Izv.
\vol 22 \yr 1984 \pages 99--172
\endref

\ref
\key N8 
\by V.V. Nikulin
\paper
Surfaces of type K3 with finite automorphism group and Picard group of
rank three
\jour Trudy Inst. Steklov
\yr 1984
\vol 165
\pages 113--142
\transl\nofrills English transl. in
\jour  Proc. Steklov Math. Inst.
\yr 1985
\vol 3
\endref

\ref
\key N9 
\by V.V. Nikulin
\paper Discrete reflection groups in Lobachevsky spaces and
algebraic surfaces
\inbook Proc. Int. Congr. Math. Berkeley 1986
\vol  1
\pages 654--669
\endref

\ref
\key N10 
\by V.V. Nikulin
\paper A lecture on Kac--Moody Lie algebras of the arithmetic type
\jour Preprint Queen's University, Canada
\vol \#1994-16,
\yr 1994 \moreref alg-geom/9412003 
\endref

\ref
\key N11 
\by V.V. Nikulin
\paper Reflection groups in Lobachevsky spaces and
the denominator identity for Lorent\-zian Kac--Moody algebras
\jour Izv. Akad. Nauk of Russia. Ser. Mat.
\vol  60
\issue 2
\yr 1996
\pages 73--106
\transl\nofrills English transl. in
\jour Izvestiya Math. \vol 60 \yr 1996 \issue 2 
\pages 305--334 
\moreref alg-geom/9503003 
\endref

\ref
\key N12 
\by V.V. Nikulin
\paper The remark on discriminants of K3 surfaces moduli as sets
of zeros of automorphic forms
\jour  Preprint \yr 1995 
\moreref alg-geom/9512018
\endref

\ref
\key N13 
\by V.V. Nikulin
\paper Basis of the diagram method for generalized reflection groups 
in Lobachevsky spaces and algebraic surfaces with nef anticanonical 
class 
\jour Intern. J. of Mathem. 
\vol 7 \issue 1 \yr 1996 \pages 71--108 
\moreref alg-geom/9405011 
\endref

\ref 
\key N14
\by V.V. Nikulin 
\paper Diagram method for 3-folds and its application to
K\"ahler cone and Picard number of Calabi-Yau 3-folds. I
\inbook Higher dimensional complex varieties:
Proc. of Intern. Confer. held in Trento, Italy, June 15-24, 1994.
\eds  M. Andreatta, Th. Peternell 
\publ de Gruyter 
\yr 1996 
\pages 261--328
\moreref alg-geom/9401010 
\endref 


\ref
\key N15 
\by V.V. Nikulin
\paper K3 surfaces with interesting groups of automorphisms
\jour  Preprint RIMS Kyoto Univ. \yr 1997 \vol 1132 
\endref

\ref
\key R
\by M.S. Raghunatan 
\book Discrete subgroups of Lie groups 
\publ Springer
\yr 1968
\endref 

\ref
\key P-SSh
\by I.I. Pjatetcki\u i-\u Sapiro and I.R. \u Safarevich
\paper A Torelli theorem for algebraic surfaces of type K3
\jour Izv. Akad. Nauk SSSR Ser. Mat.
\vol  35  \yr 1971 \pages 530--572
\transl\nofrills English transl. in
\jour Math. USSR Izv.
\vol 5 \yr 1971
\endref

\ref
\key V1 
\by \'E.B. Vinberg 
\paper On groups of unit elements of certain quadratic forms
\jour Mat. Sbornik
\yr 1972
\vol 87
\pages 18--36
\transl\nofrills English transl. in
\jour Math USSR Sbornik
\vol 16
\yr 1972
\pages 17--35
\endref

\ref
\key V2 
\by \'E.B. Vinberg 
\paper The absence of crystallographic reflection groups in Lobachevsky
spaces of large dimension
\jour Trudy Moscow. Mat. Obshch.
\vol  47 \yr 1984  \pages 68 -- 102
\transl\nofrills English transl. in
\jour Trans. Moscow Math. Soc.
\vol 47 \yr 1985
\endref

\ref
\key V3 
\by \'E.B. Vinberg 
\paper Hyperbolic reflection groups
\jour Uspekhi Mat. Nauk
\vol 40
\yr 1985
\pages 29--66
\transl\nofrills English transl. in
\jour Russian Math. Surveys
\vol 40
\yr 1985
\endref

\endRefs
\enddocument
\end